\documentclass{article}

\usepackage{amsmath, amssymb, amsthm, epsfig}
\usepackage{times}
\usepackage[noend]{algorithmic}
\usepackage{algorithm}

\author{Gregorio Malajovich\\
{\small \em Dep. de Matem\'atica Aplicada, Universidade Federal do Rio de Janeiro}\\
{\small \em Caixa Postal 68530,
Rio de Janeiro, RJ, 21945 -- BRASIL}\\
{\small \em gregorio@labma.ufrj.br,
http://www.labma.ufrj.br/\~{}gregorio}
\\
\and
Jorge P.  Zubelli\\
{\small \em IMPA}\\ {\small \em Estrada Dona Castorina 110, Rio de Janeiro, RJ, 22460-320, BRASIL}\\ {\small \em zubelli@impa.br, http://www.impa.br/\~{}zubelli}\\} 
\date{Second revision, \today}

\title{Tangent Graeffe Iteration}
\date{Revised version, August 26, 1999}

\newtheorem{theorem}      {Theorem}
\newtheorem{lemma}        {Lemma}
\newtheorem{proposition}  {Proposition}
\newtheorem{corollary}	  {Corollary}
\theoremstyle{definition}
\newtheorem{example}      {Example}
\newtheorem{definition}   {Definition}
\newtheorem{induction}    {Induction Hypothesis}

\newcommand{\condition}    {circle free \ }
\newcommand{\bydef}        {\stackrel{\text {\tiny def}}{=}}
\newcommand{\myqed}        {\begin{flushright}\qed\end{flushright}}
\newcommand{\macheps}      {\ensuremath{\epsilon_{\mathrm m}}}
\renewcommand{\Re}         { \text {Re\ }}
\renewcommand{\Im}         { \text {Im\ }}
\newcommand{\variation}[2]{ \frac{r(#2) - r(#1)}{#2 - #1} }

\newcommand{\restricted}[1]{\raisebox{-1.5ex}{$\bigg |_{#1}$} }

\newcommand{\newtext}[2]{#1}

\begin{document}
\maketitle
\begin{abstract}
  Graeffe iteration was the choice algorithm for solving univariate
  polynomials in the XIX-th and early XX-th century. In this paper,
  a new variation of Graeffe iteration is given, suitable to 
  	IEEE floating-point arithmetics of 
  modern 
  digital computers. 

	We prove that under a certain generic assumption the proposed 
	algorithm converges. We also estimate the error after $N$ iterations
	and the running cost.
\par
  The main ideas from which this algorithm is built are: classical Graeffe
  iteration and Newton Diagrams, changes of scale (renormalization), 
  and replacement of a difference technique by a differentiation one.

	The algorithm was implemented successfully and a number of numerical
	experiments are displayed.
\end{abstract}
\tableofcontents
 
\section{Introduction}

Many present day numerical algorithms have 
originated in highly ac-{\linebreak}claimed methods dating from last century
or even earlier. Such is certainly the case of Euler's method
or of Newton's method, whose numerical and theoretical consequences
still impact us today \cite{KA,BEZI,BEZII,BEZIII,BEZV,BEZIV,VAINBERG}.

Graeffe's classical method
for finding simultaneously all roots of a polynomial was
introduced independently by
Graeffe, Dandelin and Loba-{\linebreak}tchevsky~\cite{HOUSEHOLDER}.
Its simplicity,  as well as importance throughout last century 
indicate its potential as an effective numerical algorithm. 

Surprisingly, Graeffe's method has not received much attention in 
\linebreak
present day numerical computations.
Very few modern discussions about it or its applications can be found.
See the review by V.~Pan~\cite{PAN}, and also
~\cite{BP,DEDIEU,DGY,GRAU,KIRRINIS,MZ96,MZ97,NR,PAN96,PKSAH,SCHONHAGE}.

One of the main reasons for Graeffe's lack of popularity stems from 
the fact that
its traditional form leads to exponents that easily exceed the maximum
allowed by floating-point arithmetic.
Other reasons, such as the ``chaotic" behavior of the arguments
of the roots of the iterates contribute to such stigma.

\newtext{\sloppypar
Also, Graeffe iteration is a many-to-one map. It can map well-conditioned
polynomials into ill-conditioned ones, as pointed out by Wilkinson
in~\cite{WILKINSON}. We shall refer to this as `Wilkinson's
Deterioration of Condition'. 
}
{Comments on Wilkinson's deterioration of condition, etc...}

In this work we present a version of Graeffe's algorithm, 
which is well suited for floating-point arithmetic computations.
Furthermore, it has excellent complexity and memory allocation characteristics.
Our method computes both the moduli and the argument of 
all the roots, provided that certain generic conditions are satisfied.
These claims are backed by our theoretical results presented in the
next section, and proved throughout the paper, 
as well as the numerical experiments presented in the
end of the paper. 

The main ingredients in our approach are the following:

\begin{itemize}
\item The idea of renormalizing the relevant operations at each
iteration step, akin to what is done in  
dynamical systems and physics~\cite{MACKAY,MCMULLEN}.
\item The idea of using the differential of our Renormalized
 Graeffe iteration as a way of keeping the information concerning the 
argument of the roots.
\newtext{This will allow us to avoid the
harmful effects of Wilkinson's `deterioration of condition',
as discussed in Section~\ref{detcond}.}
{Wilkinson again}
\item A renormalized version of Newton's diagram that allows us to 
recognize and locate pairs of conjugate roots, as well as roots
of higher multiplicity.
\end{itemize}

The first idea mentioned above was developed in our earlier work
\cite{MZ97}, which in a certain sense laid the conceptual framework
for our present approach. It is {\em not} however essential in 
understanding the proofs presented herein. 

The second ingredient mentioned above is explained and
motivated in Section~\ref{tangra}.
In rather vague terms it could be compared to the advantage of 
using derivatives, when those are available,
as compared to using differences. 
\newtext{This idea can be traced back to Brodetsky and
Smeal~\cite{BS} in 1924, in a more ad-hoc fashion. We are not
aware of recent applications of that method in modern literature.}
{New reference}

Finally, the concept of Newton's diagram, as well as the power of 
Graeffe's method was present throughout Ostrowski's 
masterpiece~\cite{OSTROWSKI}. 
While writing the present paper we could not help but wonder what 
would have been the outcome of that research if he had available at that
point the present day technology of high speed computers.

\newtext{We wish to thank two anonymous referees for their
comments and for suggesting some extra references such as ~\cite{BS},
which we were not aware of in the first draft.}{Thanks}

\subsection{Main Result}

	We will introduce an algorithm for solving real and complex
        univariate polynomials.
	The following genericity condition will be required at input:

\begin{definition}
   A {\em real} polynomial $f$ will be called {\em \condition} if, and only if, 
   for any couple $\zeta$, $\xi$
   of distinct roots of $f$, one has either $|\zeta| \ne |\xi|$, or $\zeta = \bar \xi$.
\end{definition}
\begin{definition}
   A {\em complex} polynomial $f$ will be called {\em \condition} if, and only if, for 
   any couple $\zeta$, $\xi$
   of distinct roots of $f$ one has  $|\zeta| \ne |\xi|$.
\end{definition}

   It is obvious that given any real polynomial $f$, one can obtain a \condition
   polynomial by (pre)composing it with a conformal transform of the form:
\[
   \begin{array}{crcl} \varphi: & \bar \mathbb C & \rightarrow & \bar \mathbb C \\
        & x & \mapsto & \varphi(x) =  \frac{x \cos \theta - \sin \theta}{x \sin \theta + \cos \theta}
   \end{array}  
\]
   and then clearing denominators; for all but a finite number of $\theta \in (-\pi, \pi]$, the
   resulting polynomial
\[
   \tilde f(x) = (x \sin \theta + \cos \theta)^d 
     f \left(
     \frac{x \cos \theta - \sin \theta}{x \sin \theta + \cos \theta} 
     \right)
\] 
   is \condition\negthickspace\negthickspace\newtext
   {, where $d$ is the degree of $f$}
   {$d$ was not defined.}.
\medskip
\par
   Tangent Graeffe Iteration will be shown to converge for all \condition polynomials; 
   Given an arbitrary polynomial, one can first find all zero roots (in the obvious way), 
   apply a random conformal transform, then 
   Tangent Graeffe Iterations, and finally recover the roots of the original polynomial.
\medskip
\par

When counted with  multiplicity, the roots of a \condition polynomial
can be canonically ordered by:
\begin{enumerate}
\item $| \zeta_1 | \le \cdots \le | \zeta_d | $ 
\item In the real case, 
\subitem 2.1. \ If $|\zeta_i|=|\zeta_{i+1}|$ then $\zeta_i=\overline{\zeta_{i+1}}$.
\subitem 2.2. \ If $i=1$ or $|\zeta_{i-1}|< |\zeta_{i}|$ then $\Im{\zeta_{i}} \ge 0$.

\end{enumerate}
\medskip
\par
   If we assume that all the arithmetical operations are performed
exactly (including transcendental), the mathematical properties of the 
algorithm can be summarized by:

\begin{theorem}\label{main}
        Let $f$ be a real (resp. complex) \condition degree $d$ polynomial,
        not vanishing at 0. Denote by 
$\zeta$  the  vector of all the roots of $f$ with multiplicity  
        canonically ordered as above.

        Then, a total of $N$ iterations of Renormalized Tangent Graeffe 
(Algorithm~\ref{alg:solve})
	produces  $\zeta^{(N)} \in \mathbb C^d$, 
such that
\[ \zeta^{(N)} \longrightarrow \zeta \ \ \ \mbox{  (as $ N \rightarrow \infty$).  } \]
\par
	The running time for each iteration is $O(d^2)$ exact arithmetic operations
	(including transcendental operations). The relative 
	\newtext{truncation}{for clarity}
	error bound in
        each coordinate after $N$ iterations
	is $2^{-2^{N-C}}$, where $C$ depends on $f$.
\par
\end{theorem}

\subsection{What the  Graeffe Iteration is; Its historical weaknesses}

In this section we shall briefly review the main ideas behind the
method and describe also some of its weaknesses.

	Graeffe iteration maps a degree $d$ polynomial $f(x)$ into the
degree $d$ polynomial 
\[
Gf(x) = (-1)^d f(\sqrt{x}) f(-\sqrt{x}) \mbox{ .} 
\] 
If $\zeta_1, \zeta_2, \dots \zeta_d$ are the roots of $f$, then the
roots of $Gf$ are $\zeta_1 ^2, \zeta_2 ^2, \dots \zeta_d ^2$ 
\medskip
\par
	Assume that $g = G^N f$
is the $N$-th iterate of $f$. Then, assuming
that $f$ is monic, the 
coefficients of $g(x)=g_0 + g_1 x + \cdots +  g_d x^d$ satisfy: 
\[
\begin{array}{lcl}
g_0  & = & (-1)^d \ \sigma_d \ \left( \zeta_1 ^{2^N}, \zeta_2 ^{2^N}, \cdots,  \zeta_d ^{2^N} \right) \\
g_1  & = & (-1)^{d-1}\  \sigma_{d-1} \ \left( \zeta_1 ^{2^N}, \zeta_2 ^{2^N}, \cdots,  \zeta_d ^{2^N} \right) \\
& \vdots & \\
g_j  & = & (-1)^{d-j}\  \sigma_{d-j} \ \left( \zeta_1 ^{2^N}, \zeta_2 ^{2^N}, \cdots,  \zeta_d ^{2^N} \right) \\
& \vdots & \\
g_d  & = & \sigma_0 = 1 
\end{array}
\]
where $\sigma_j$ is the $j$-th elementary symmetric function.
In the particular case that $|\zeta_1| < |\zeta_2| < \dots < | \zeta_d|$, we can further approximate
\[
\begin{array}{lcl}
g_0  & = & (-1)^d \zeta_1 ^{2^N} \zeta_2 ^{2^N} \dots \zeta_d ^{2^N}  \\
g_1  & \simeq & (-1)^{d-1} \zeta_2 ^{2^N} \dots \zeta_d ^{2^N} \\
& \vdots & \\
g_j  & \simeq & (-1)^{d-j} \zeta_{j+1} ^{2^N} \dots \zeta_d ^{2^N} \\
& \vdots & \\
g_d  & = & \sigma_0 = 1 
\end{array}
\]
Hence, it is possible to determine 
\[ \zeta_j^{2^N} \simeq -\frac{g_j}{g_{j+1}} \mbox{ .} \] 

\medskip
\par
	We stress two main weaknesses. 
	The first big weakness of classical Graeffe iteration is coefficient
growth. As the coefficients of $g_j$ grow doubly exponentially in the number
of iterations, the exponent (not the mantissa) of the floating-point 
system gets overflowed:
 
\begin{example} \label{overflow1}
  Let $f$ have roots $1$, $2$, $3$, $4$. Then the $N$-th Graeffe iterate of $f$ has roots
  $1$, $2^{2^N}$, $3^{2^N}$, $4^{2^N}$. The coefficient $g_0$ is $24^{2^N}$. If $N=8$, 
  then $g_0$ is approximately $1.68 \times 2^{1173}$, while IEEE double precision numbers 
  used in most modern computers cannot contain floating point values more than $2^{1024}$,
  since the exponent is represented by 11 bits (sign included) \cite{HIGHAM}. (As a matter
  of fact, the representation is a little more complicated, as it allows
 for `subnormal' numbers \cite{DEMMEL}). Therefore,
  we would have an overflow when computing the $8$-th Graeffe iterate of $f$.
\end{example}
\begin{example} \label{overflow2}
  On the example above, assume that $f$ would have an additional root $1.01$. Namely,
\[
  \begin{array}{rcl}
  f(x) &=& (x-1)(x-1.01)(x-2)(x-3)(x-4) \\
       &=& -24.24 + 74.5x -85.35 x^2 + 45.1 x^3 -11.01 x^4 + x^5
  \end{array}
\]
   We will show that 8 Graeffe iterations are not enough to compute
   the first root (namely 1) to an accuracy of $10^{-4}$. However, as 
   shown in Example~\ref{overflow1}, 8 iterates are enough to 
   overflow the IEEE double precision number system.
\par
   Indeed, the first root is obtained as:
\[
   \begin{array}{rcl}
   \zeta ^{2^N} \simeq - \frac{g_0}{g_1} &=& 
   - \frac{1.01 ^{2^N} 24^{2^N}}{(1^{2^N} +1.01 ^{2^N}) 24^{2^N}} \\
   &=& \frac{1}{1 + 1.01^{-2^N}} \\
   &\simeq& 0.927  \ \ \  \mbox{( for $N=8$).  }   	
   \end{array}
\]
   Thus,
\[
   \zeta \simeq 1 - 2.9 \times 10^{-4}  
\] 
\par
   The error obtained is therefore larger than $10^{-4}$.
\end{example}
\medskip
\par
	The introduction of the idea of renormalization allows us to
avoid coefficient growth, and replace a diverging algorithm by a 
convergent one (See Section~\ref{renorm} and also ~\cite{MZ97}). 
Alternative approaches for the number range growth are suggested in
~\cite{HENRICI} and in~\cite{GRAU}.
\par
A certain geometrical invariant of the polynomial, 
the  {\em limiting  Newton diagram}, appears naturally in the context of 
Renormalized Graeffe Iteration. It allows to recover the information
about multiple roots and pairs of roots.  (See ~\cite{OSTROWSKI}). 
\medskip
\par
	In Section~\ref{newton}, we give a procedure to obtain
the limiting Newton diagram of a given polynomial. It is effective in
the sense that, if we can bound the separation 
\[ \max_{|\zeta_i| > |\zeta_j|} \frac{|\zeta_i|}{|\zeta_j|} \mbox{ , } \]
then we can effectively identify the
multiple roots and pairs of roots. It will converge, and eventually provide
the list of multiple roots and pairs for any circle-free polynomial,
in a finite (but unknown, not effective) number of iterations.

\medskip
\par
	The other big weakness of classical Graeffe iteration is the
	fact that it returns the moduli of the roots, but not the actual
	roots. As a matter of fact, information about the argument of
	the roots is lost, and should be recovered by other means:
\begin{example}
	Consider polynomials
	$f(x) = x^2 - 2x + 1$, $g(x) = x^2 - 1$, $h(x) = x^2 + 1$.
	After two Graeffe iterations, all the three polynomials are
	mapped into $f(x)$.
\end{example} 

	Many algorithms have been proposed to recover the arguments
~\cite{PAN}.
	In this paper, we will differentiate the Graeffe iteration operator, 
	and obtain
	an iteration defined on the appropriate tangent bundle. This
	new operator will define a mapping between $1$-jets of polynomials.
        By the latter we mean expressions of the form 
$f(x) + \epsilon \dot f(x)$, where
	$\epsilon$ is a formal parameter. This procedure is discussed in
	section~\ref{tangra}. \newtext{
	In the end of the same section, we shall discuss the stability
	properties of this process.}{Wilkinson again !}
\par
	In section~\ref{numerics}, we compare the numerical behavior of
	Renormalized Tangent Graeffe Iteration to other publicly 
	available algorithms. 
\section{Renormalizing Graeffe} 
\label{renorm}

\subsection{The Renormalized Graeffe Iteration}
\medskip
\par
	Example~\ref{overflow2} shows a typical behavior of 
        classical Graeffe iteration performed by
	digital computers \cite{HENRICI}. In order to avoid that sort of overflow,
        the authors introduced in \cite{MZ97} the {\em Renormalized Graeffe Iteration}.
	Although the details and the mathematical foundations of the algorithm are
	described in \cite{MZ97}, to keep the present work self-contained,
	we give below a very short description of the main
	ideas:
\medskip
\par
	One should consider the computation of $g = G^N f$ as 
	divided in several {\em renormalization levels}.
\par
\centerline{
\begin{tabular}{lll}
Level 0  &  Coefficients of $f$ \\
Level 1  &  Coefficients of $Gf$ \\
Level 2  &  Coefficients of $G^2f$ \\
         & \multicolumn{1}{c} \vdots \\
Level N &  Coefficients of $G^Nf$ \\
\end{tabular}} 
\medskip
\par
    At renormalization level $N$, all coefficients $g_j$ of $g=G^N f$ 
    should be represented in
    coordinates 
 \[ r_j^{(N)} = - 2^{-N} \log |g_j| \mbox{ ,} \] 
and 
    \[ \alpha_j^{(N)} = g_j / |g_j| \mbox{ .} \] 
    Therefore, we shall obtain convergence
    of the radial coordinates $r_j^{(N)}$ (at least in the case of roots 
    of different moduli). The dynamics
    of the angular coordinates $\alpha_j^{(N)}$ is typically chaotic.
\par
	In order to pass from level $N$ to level $N+1$, a {\em Renormalized Graeffe Operator}
    was defined in~\cite{MZ97}. Intermediate computations were performed in coordinates
    $r_j^{(N)}$ and $\alpha_j^{(N)}$ by means of renormalized arithmetic operations. For instance, the 
    renormalized sum $(r,\alpha)$ of $(r_1,\alpha_1)$ 
and $(r_2 , \alpha_2)$ can be defined 
    (in renormalization level $N$) by
\[
    r =  - 2^{-N} \log | \alpha_1 e^{- 2^N r_1}
                    + \alpha_2 e^{- 2^N r_2} |   
\]
\[
    \alpha = \frac{ 
                    \alpha_1 e^{ - 2^N r_1}
                    + \alpha_2 e^{ - 2^N r_2}    
              }
              {|  
                    \alpha_1 e^{ - 2^N r_1}
                    + \alpha_2 e^{ - 2^N r_2}    
             | }
\]
 
\par
	Renormalized sum can be computed without
   overflow by the formula in Algorithm~\ref{alg:rensum}.
This is a simplified, non-optimal version of 
renormalized sum. Notice that one or two of the inputs can be the renormalization
of $0$, i.e., $\infty$. Under the usual conventions, $\infty$ is greater than any
real number. Therefore, if only one of the arguments is $\infty$, the correct result
will be returned.

\begin{algorithm}[!ht]
\caption{RenSum ( $r_1$, $\alpha_1$, $r_2$, $\alpha_2$, $p$ )}
\label{alg:rensum}
\begin{algorithmic}
\STATE
\COMMENT {~ It is assumed that $r_1$ and $r_2$ are real numbers or $+ \infty$, and that
           $|\alpha_1| = |\alpha_2| = 1$. The number $p$ should be equal
	   to $2^N$, where $N$ is the renormalization level. This routine
           computes (in renormalized coordinates !) the sum of 
	   $\alpha_1 e^{-p r_1}$ and $\alpha_2 e^{-p r_2}$.~}
\IF {$r_1 = r_2 = + \infty$} \STATE \textbf{return} $+\infty, 1$ \ENDIF 
\STATE $\Delta \leftarrow r_2 - r_1$
\IF {$\Delta \ge 0$}
    \STATE $t \leftarrow \alpha_1 + \alpha_2 e^{ - p \Delta }$
    \STATE \textbf{return} $r_1 - \frac {\log(|t|) }{ p}, t/|t|$
    \ELSE
    \STATE $t \leftarrow \alpha_2 - \alpha_1 e ^{p \Delta}$
    \STATE \textbf{return} $r_2 - \frac{\log(|t|) }{ p}, t/|t|$
    \ENDIF
\end{algorithmic}
\end{algorithm}

\medskip
\par
	A few extra mathematical ideas related to the  
	renormalized Graeffe operators, as well as 
        some other mathematical results can be found in \cite{MZ97}.  
 
\subsection{The Renormalized Newton Diagram}

	The first goal of this section is to introduce the concept of
   	Renormalized Newton diagram, which is going to play a key 
	role in the practical implementation of the algorithm 
	discussed in this paper. The second is to prove a convergence
	result based on such idea using some earlier results of Ostrowski's.
	
	We start by reviewing the concept of Newton diagram, which
	has been used extensively by Ostrowski, Puiseux and Dumas, 
	among others.

	Let 
		\[ f = \sum_{i=0}^{d} f_{i} x^i  \mbox{ ,} \]
	be a degree $d$ polynomial.  As before, we denote by $g=G^N f$ the
	$N$-th Graeffe iterate of $f$.

	We order the roots of $f$ in nondecreasing order of their moduli,
	to wit:
	\begin{equation} \label{ineqs}
	 | \zeta_{1} |  \le \dots \le  | \zeta_{d} | \mbox{ .} 
	\end{equation}

	If the above inequalities are all strict, then 

        \[ \lim_{N\rightarrow \infty}
                2^{-N} \log\left(\frac{|g_{i}|}{|g_{i+1}|}\right) =
                \log(|\zeta_{i+1}|) \mbox{ .} \]

	For each $N$, consider the piecewise linear function 
	$r^{(N)}:[0,d] \rightarrow {\mathbb{R}}\cup \{+\infty\} $
	satisfying 
		\[ r^{(N)}(i)\bydef -2^{-N}\log | g_i| \mbox{ .} \]
\begin{figure}
\centerline{\epsfig{file=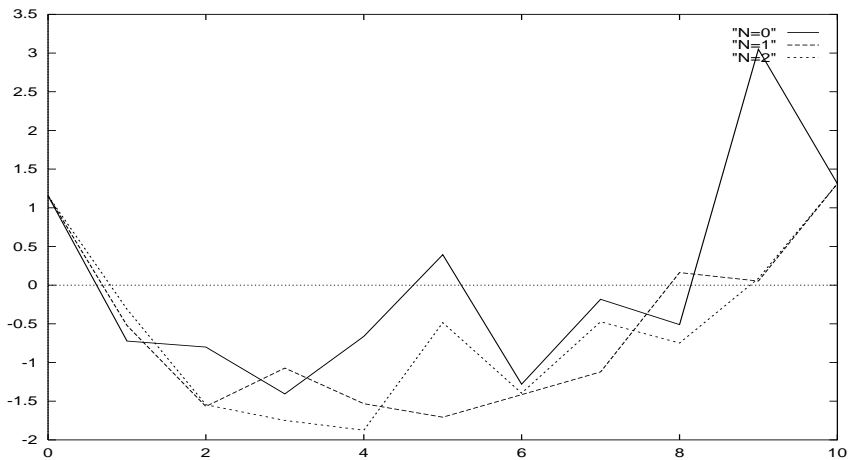, height=6cm, width=12cm}}
\caption[ ]{\label{fig3}The function $r^{(N)}(i)$, for $N = 0, 1, 2$}
\end{figure}
	Notice that under the above assumptions, $r^{(N)}$ is convex for 
	sufficiently large $N$. (See figure~\ref{fig3}).
	Indeed, since $|\zeta_{i+1}|>|\zeta_{i}|$, 
        the inclinations satisfy (for large $N$)
	\[ r^{(N)}(i+1) - r^{(N)}(i) \ \ge \ r^{(N)}(i) -  r^{(N)}(i-1) \mbox{ .} \]

	It is easy to see that if two consecutive roots, say $\zeta_{i}$ 
	and $\zeta_{i+1}$, have approximately the 
	same absolute value, then the three corresponding points 
	\[ \left(i-1,r^{(N)}(i-1)\right), 
	\left(i,r^{(N)}(i)\right), 
	\left(i+1, r^{(N)}(i+1)\right) \mbox{ } \]
	will be approximately aligned.
	Furthermore, the functions $r^{(N)}$ converge to a piecewise linear
	convex function. 

\medskip
\par
	However, if the inequalities in (\ref{ineqs}) are not strict, the
        functions $r^{(N)}$ may fail to converge.
\begin{example}
	Let $f(x) = (x-1)(x-e^{i \theta})$. Then its $N$-th Graeffe
        iterate is
\[
	g(x) = x^{2} - (1 + e^{2^N i \theta}) x + e^{2^{N} i \theta}
\]
	Therefore, we have $r^{(N)} (0) = r^{(N)}(2) = 0$, but we
        also have 
\[
r^{(N)}(1) = -2^{-N} \log|1 + e^{2^N i \theta}|
        = -2^{-N} \log |2 \cos 2^{N-1}\theta|\mbox{ \ .}
\]
	Depending on the choice of $\theta$, this last value can range
        anywhere from $-2^{-N} \log 2$ to $+ \infty$.
\end{example} 
\medskip
\par
\begin{figure}
\centerline{\epsfig{file=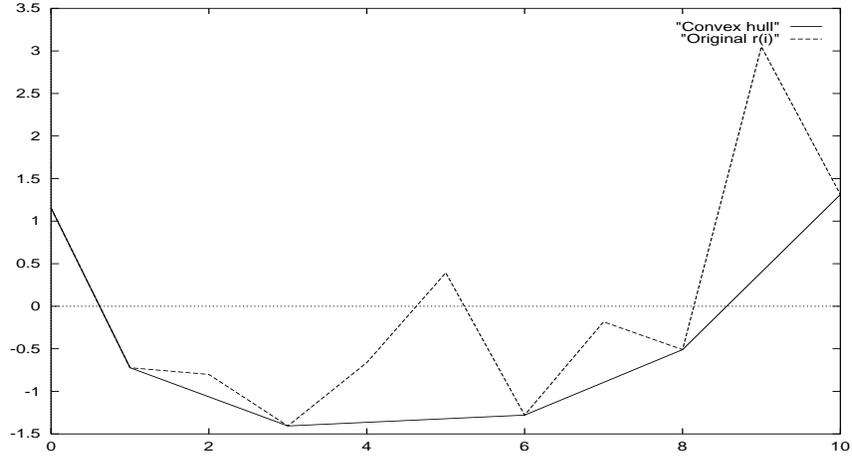, height=6cm, width=12cm}}
\caption[ ]{\label{fig4}The function $r^{(N)}(i)$ and its Convex Hull}
\end{figure}

        This is one of the reasons for introducing the 
        convex hull of each $r^{(N)}$, that will 
	be subsequently called the {\em Renormalized Newton Diagram}. 
	(See figure~\ref{fig4}). 
\par
        Our approach has the
	advantage of simplifying some of the arguments 
        by Ostrowski in ~\cite{OSTROWSKI} by providing
        plain convergence of Renormalized Newton Diagrams; 
        however, we will quote several of the results
        by Ostrowski in the sequel.	
\medskip
\par
	One of the major goals of Ostrowski in ~\cite{OSTROWSKI}
        was to obtain effective bounds for the moduli of roots.
        This was possible by introducing of the {\em majorant}
        of a given polynomial:
\begin{definition}
	A majorant of a given polynomial is any other 
	polynomial, of same degree, with nonnegative coefficients greater 
	than or equal to the given polynomial's coefficients.
\end{definition}
	The first step of Ostrowski's construction is 
	Newton's majorant:
\begin{definition}
	A polynomial $A = \sum_{i=0}^{d} A_{i} x^i$  with nonnegative
	coefficients is called {\em normal} if the following conditions
        hold:
	\begin{enumerate}
	\item If $A_{i}>0$ and $A_{j}>0$ for $i>j$, then
		$A_{l} > 0$ for all $i < l < j$. 
	\item For $l=1,\dots,d-1$, 
	\[ A_{l}^2 \ge A_{l-1} A_{l+1}\mbox{ .} \]
	\end{enumerate}
\end{definition}
	A normal majorant $T = \sum T_i x^i$ of $f$ 
	is called {\em minimal} if for any 
	other majorant $T'$ of $f$ we have
\begin{equation} \label{star1}
 T_{j} \le T'_{j} \mbox{ , $j=0, 1, \dots, d$.} 
\end{equation}
Notice that Condition~2 above
means that the graph of the points of the form $(l,-\log(r_{l}))$, for $l=0, \dots, d$
is convex.

	In the language of majorants,

\begin{proposition}[Ostrowski\cite{OSTROWSKI}]
Any polynomial 
\[ f = \sum_{j=0}^{d} f_{j} x^j \]
possesses a unique minimal normal majorant, 
\[ {\cal M}_{f} = \sum_{j=0}^{d} T_{j} x^{j} \mbox{ .} \]
\end{proposition}
\newtext{
\medskip
\par
The polynomial ${\cal M}_{f}$ will be called {\em Newton Majorant}
of the polynomial $f$.}{Definition of Newton Majorant}

\medskip
\par
The result can be  
proved by using the convex hull $\phi$ of the function $-\log |f_i|$.
and constructing the polynomial
$ {\cal M}_f$ as the polynomial with positive coefficients 
$({\cal M}_{f})_j = e^{-\phi(j)}$.
 We refer the
interested reader to Ostrowski's work \cite{OSTROWSKI,OSTROWSKI2}.
\medskip
\par

We remark that if the polynomial has roots of strictly increasing moduli,
then the coefficients $T_{i}$ of the Newton's majorant of $g=G^N f$
coincide with $|g_{i}|$ for $N$ sufficiently large and $ i=0,1, \dots, d$.
\medskip
\par 
However, the introduction of the Newton Diagram allows us to consider
the general situation of possibly many roots
of same moduli. As before, we order them in nondecreasing order 
and consider the indices
\[ i_{0}=0 <  i_{1} < i_{2} < \dots < i_{l} < i_{l+1}=d \]
as $0$, $d$ and exactly those integers $i$ between $1$ and $d-1$ such that
$|\zeta_{i}| < |\zeta_{i+1}|$. This way, we have that 
\[ 
|\zeta_{i_{j-1}}| < |\zeta_{i_{j-1}+1}| = \dots =
|\zeta_{i_{j}}| < |\zeta_{i_{j}+1}| \mbox{ .} 
\]

The fundamental result, in this case is 
\[ 
\lim_{N\rightarrow \infty}
2^{-N} \log  \frac{|g_{i_{j}}|}{|g_{i_{j+1}}|} 
= ( i_{j+1}-i_{j}) \log | \zeta_{i} | \mbox{ ,} \]
for $i_{j} < i \le i_{j+1}$ and $0\le j \le l$.
(c.f. equation (79.8) of \cite{OSTROWSKI2}). In the language
of Renormalized Newton Diagrams, that very same equation
can be written as:
\[
\log |\zeta_{i}| = \lim_{N\rightarrow \infty}
\frac {r^{(N)} (i_{j+1}) - r^{(N)} (i_j)}{i_{j+1} -  i_j} 
\mbox{ ,} \ \ \ \ \mbox{ for  $i_{j} < i \le i_{j+1}$.} 
\] 
\medskip
\par

As remarked by Ostrowski, the above formulae are only useful in the
determination of the moduli of the roots if we know ``a priori" the
values $i_{1} < i_{2} < \dots < i_{l}$.  This is obviously not
the case in most applications. Instead, Ostrowski's results
provide effective bounds for the convergence of the $r^{(N)}$,
and thus for the values of the $|\zeta_i|$' s.  
\medskip
\par

\begin{theorem}[Ostrowski\cite{OSTROWSKI},Theorem IX.3]
Let  \[ \varrho(\nu) \bydef 1 - 2^{-1/\nu} \mbox{ ,} \]
and 
\[ R_{\nu}^{(N)} \bydef \frac{T_{\nu-1}^{(N)}}{T_{\nu}^{(N)}} \mbox{ ,} \]
then 
\begin{equation}
\varrho(\nu) < \frac{|\zeta_{\nu}|^{2^N}}{R_{\nu}^{(N)}} < 
\frac{1}{\varrho(d-\nu+1)} \mbox{ \ \ \, \ \ \ $\nu =1, \dots, d$.} 
\end{equation}
\end{theorem}

As a consequence of the above estimate, Ostrowski gets the following
bound 
\[  (2d)^{-2^{-N}} < \frac{|\zeta_{\nu}|}{(R_{\nu}^{(N)})^{2^{-N}}} < 
(2d)^{2^{-N}} \]

\begin{corollary}
If $r^{(N)}(i)$ denotes the $i$-th ordinate of the $N$-th Renormalized
Newton Diagram, then
\[ \lim_{N\rightarrow \infty} \left(r^{(N)}(i) - r^{(N)}(i-1)\right) = 
\log |\zeta_{i}|  \]
for $i=1,\dots,d$. 
Furthermore, the error is bounded from above by  
\[ 2^{-N} \log(2d) \mbox{ .}  \]
\end{corollary}
\par
This is indeed a strong result, since nothing is assumed on 
the coefficients or the roots of the original polynomial. However, it is
possible to get a better error bound by assuming 
a minimal separation on the moduli:
\[
\min_{ |\zeta_i | > |\zeta_j| } \frac {|\zeta_i|}{|\zeta_j|} > 1 + \epsilon
\] 
for some $\epsilon >0$.
\par
	We note that in the above formula, if $\epsilon$ is well defined
(i.e. there are at least two roots of different modulus) then 
it is non-zero.

%
%
%
%
%
%
%
\section{Computation of the Newton Diagram}
\label{newton}
\subsection{Algorithm and Main Statements}
	The main issue in this section is the following:

	We are given a certain polynomial $g$, obtained 
	after a few Graeffe iterations of a polynomial $f$.
	The roots of $g$ are $Z_1$, \dots , $Z_d$, and we
	order them so that:
\[
	|Z_1| \le |Z_2| \le \dots \le |Z_d|
\]
\par
	We want to know which of the inequalities are strict.
	We do not know the actual value of the $Z_i$'s, we know
	only the coefficients of $g$, i.e., the symmetric 
	functions of the $Z_i$'s.
\par
	We are also ready to assume that
\[
	R \bydef \min_{|Z_{i+1}| > |Z_i|} \frac{|Z_{i+1}|}{|Z_i|}
\]
	is a large real number. Indeed, if $\zeta_1$, \dots ,
	$\zeta_d$ are the roots of $f$, always ordered such as
	$|\zeta_1| \le \dots \le |\zeta_d|$, then 
\[
	\rho \bydef \min_{|\zeta_{i+1}| > |\zeta_i|} \frac{|\zeta_{i+1}|}{|\zeta_i|}
\]
	is always strictly greater than one. Hence, given any $A > 0$,
        by performing $N \ge \log_2 \frac{\log A}{\log \rho}$  
	iterations, we can assume that $R = \rho^{2^N} \ge A$. 
\medskip
\par
	Recall that the Renormalized Newton Diagram of $g$ is the
	convex hull of the function $i \mapsto -2^{-N} \log g_i$.
	As $N$
	grows, the Renormalized Newton Diagram of $g$ converges 
	to the convex hull of 
\[ i \mapsto -\log |f_0| + \sum_{j \le i} \log |\zeta_i| \mbox{ .} \]
	However, we want to be able to decide in {\em finite time} what are the
	sharp corners of the convex hull of 
        $i \mapsto \log |f_d| + \sum_{j \le i} \log |\zeta_j|)$. As in the
	preceding section, we write:
\[
	r_i = -2^{-N} \log |g_i|
\]
	where $g$ is the $N$-th iterate of $f$. Notice that we dropped
	the superscript $N$ of $r_i^{(N)}$.
 
\begin{proposition} \label{prop:cvxhull}
	Let $g = G^N f$, where $f$ is a degree $d$ polynomial and $G$ denotes the 
	Graeffe iteration. Let $Z$, $\zeta$ and $\rho$ be as above. Assume
	that  
\[ N > 3 + \log_2 \frac{ d \log 2 } {\log \rho} \mbox{ .} \]
	Then, Algorithm~\ref{alg:cvxhull} below with input 
        $N$, $d$, $r$, $\rho$ produces the list of the sharp
	edges of the convex hull of
        $i \mapsto \sum_{j \le i} \log |\zeta_i|)$.
\end{proposition}

	Remark: Algorithm ~\ref{alg:cvxhull} has running time $O(d)$.

\begin{algorithm}[!ht] 
\caption{Strict Convex Hull ( $N, d, r, \rho$ )}
\label{alg:cvxhull}
\begin{algorithmic}
\STATE
\STATE
\COMMENT {~ Create a list $\Lambda$, containing initially the element 
           $\Lambda_0=0$}
\STATE  $j \leftarrow 0$;
\STATE  $\Lambda_{j} \leftarrow 0$ ; 
\STATE
\STATE
\COMMENT {~The error bound below will follow from Lemma~\ref{cvx4}.} 
\STATE  $R \leftarrow \rho ^ {2^N}$
\STATE	$E \leftarrow \frac{1}{2} \left(2^{-N+1} \log (2^d + 2^d R^{-1})
           - 2^{-N+2} \log(1 - 2^d R^{-1}) + \frac{\log{\rho}}{2}
	  \right)$ 
\STATE
\STATE
\COMMENT {~Now, we will try to add more points to the list $\Lambda$. At
	  each step, we want to ensure that we have always a convex set.}
\STATE
\FOR { $i \leftarrow 1$ to $d$ } 
\STATE
\STATE 
\COMMENT {~We discard all the points in $\Lambda$ 
          that are external to the convex hull of $\Lambda$ and the 
          new point.
          Let $\Lambda_j$ be the last element of $\Lambda$}
\STATE
\STATE
\WHILE   { $j>0$ and   
           $\frac{ r_{\Lambda_j} -r_{\Lambda_{j-1}} } 
                 {\Lambda_j - \Lambda_{j-1}}
            >
           \frac{ r_{i} -r_{\Lambda_j} }
                 {i - \Lambda_{j}}
	   - E$}
\STATE    $j \leftarrow j-1$   
\ENDWHILE
\STATE
\STATE
\COMMENT {~Now, we append the point $i$}
\STATE  $j \leftarrow j+1$ 
\STATE  $\Lambda_{j} \leftarrow i$; 
\ENDFOR
\STATE
\STATE {\bf  Return} $(\Lambda_0, \cdots, \Lambda_j)$ 

\end{algorithmic}

\end{algorithm}

\medskip
\par
\subsection{Some Estimates about Symmetric Functions}
	In order to prove Proposition~\ref{prop:cvxhull}, we need
	a few estimates about symmetric functions. First of all,
	let $I = \{ i : |Z_i| < |Z_{i+1}| \} \cup \{ 0, d \} $ be the set of sharp
	corners of the limiting Renormalized Newton Diagram. As before,
	let $\sigma_k$ denote the $k$-th elementary symmetric function,
\[
	\sigma_k(Z) = \sum_{j_1 < \dots < j_k} Z_{j_1} \dots Z_{j_k}
\]
	Then, 
\begin{lemma} \label{cvx1}
	For $i \in I$ we have
\[
	\sigma_{d-i} (Z) = Z_{i+1} Z_{i+2} \dots Z_d (1+c) \mbox{ ,} 
\]
	where $|c| \le \left( \binom{d}{i} - 1 \right) R^{-1}
	\le 2^d R^{-1}$	
\end{lemma}

\begin{proof}[Proof of Lemma~\ref{cvx1}:]
	Write 
\[
	\sigma_{d-i} (Z) = \sum _{ j_1 < \dots < j_{d-i} } Z_{j_1} \dots Z_{j_{d-i}}
\]
\par
	In the sum above, $|Z_{j_1} \dots Z_{j_{d-i}}|$ $\le$ 
	$R^{-1} |Z_{i+1} Z_{i+2} \dots Z_d|$ 
        for any \linebreak choice of $j_1$, \dots $j_{d-i}$
	except $i+1, \dots d$. Since there are $\binom{d}{i} - 1$ other terms,
	we obtain that:
\[
        |\sigma_{d-i} (Z) - Z_{i+1} Z_{i+2} \dots Z_d |
	\le 
        \left( \binom{d}{i} - 1 \right) R^{-1}
        |Z_{i+1} Z_{i+2} \dots Z_d|
\] 	
\end{proof}

\begin{definition}
	We will say that $i_1$ and $i_2$ are {\em successive elements} of $I$
	if and only if:
\begin{enumerate}
	\item $i_1 \in I$
	\item $i_1 < i < i_2 \Rightarrow i \not \in I$
	\item $i_2 \in I$
\end{enumerate}
\end{definition}

\begin{lemma} \label{cvx2}
	Let $i_1$ and $i_2$ be successive elements of $I$, and let 
	$i_1 < l < i_2$. Then
\[
	| \sigma_{d-l}(Z) | 
	\le
	\left(
	\binom{i_2 - i_1}{i_2 - l} + c'
	\right)
	|Z_{i_2}|^{i_2-l} |Z_{i_2 + 1}| \dots |Z_d|
\]	
	with $c' \le \left( \binom{d}{i} - \binom{i_2 - i_1}{i_2 - l} \right) 
              R^{-1} < 2^d R^{-1}$
\end{lemma}

\begin{proof}[Proof of Lemma~\ref{cvx2}:]
	Write
\begin{eqnarray*}
	\sigma_{d-l}(Z) &=& \sum_{j_1 < \dots < j_{d-l}} Z_{j_1} \dots Z_{j_{d-l}} \\
	                &=& {\sum}' Z_{j_1} \dots Z_{j_{d-l}}
			   + {\sum}'' Z_{j_1} \dots Z_{j_{d-l}}
\end{eqnarray*}
	where $\sum'$ ranges over the $j$ such that $i_1 < j_1 < \dots < 
        j_{i_2 - l }< i_2+1$
	and $j_{i_2-l+1} = i_2+1$, \dots , $j_{d-l} = d$. Of course,
	$\sum''$ ranges over all the other terms.
\par
	We can rewrite $\sum'$ as:
\[
	{\sum}' = Z_{i_2+1} Z_{i_2 + 2} \dots Z_d
	        \left( \sum_{i_1 < j_1 < \dots < j_{i_2 - l} \le i_2} 
                     Z_{j_1} \dots Z_{j_{i_2-l}} \right)
\]
	Hence,
\[
	|{\sum}'| \le \binom{i_2 - i_1}{i_2 - l} 
                    |Z_{i_2}|^{i_2-l} |Z_{i_2 + 1}| \dots |Z_d|
\]
\par
	The terms in $\sum''$ are all smaller than 
	$R^{-1} |Z_{i_2}|^{i_2-l} |Z_{i_2 + 1}| \dots |Z_d|$. \linebreak
        Since there
	are $ \binom{d}{i} - \binom{i_2 - i_1}{i_2 - l}$ of them,
\[
	|{\sum}''| < \left( \binom{d}{i} - \binom{i_2 - i_1}{i_2 - l} \right)
              R^{-1}|Z_{i_2}|^{i_2-l} |Z_{i_2 + 1}| \dots |Z_d|
\]
	Adding those two bounds, we obtain indeed:	
{\small
\[
        | \sigma_{d-l}(Z) |
        \le
        \left(
        \binom{i_2 - i_1}{i_2 - l} + 
              \left( \binom{d}{i} - \binom{i_2 - i_1}{i_2 - l} \right)
              R^{-1}
        \right)
        |Z_{i_2}|^{i_2-l} |Z_{i_2 + 1}| \dots |Z_d|
\]
}
\end{proof}	

	The estimates above can be converted into `logscale' estimates:
\begin{lemma}{\label{cvx3}}
	Let $i_1$, $i_2$ be successive elements of $I$, and let 
	$i_1 < l < i_2$. Then the following three equations are true:
\begin{enumerate}
\item
\[
	\variation{i_1}{i_2}
	= \log | \zeta_{i_2} | + 2^{-N+1} \log(1+c)
\]
	with $|c| \le \left(
              \max \left( \binom{d}{i_1} , \binom{d}{i_2} \right) 
	- 1 \right)R^{-1} < 2^d R^{-1}$. 
\item
{ \small
\[
	\variation{l}{i_2}
	\le
	\log |\zeta_{i_2}| + 2^{-N} \log \left( \binom{i_2-i_1}{i_2 -l} + c'
	\right) + 2^{-N} \log |1+c| 
\]}
	where $|c| \le \left( \binom{d}{i_2} - 1 \right) R^{-1}
	\le 2^d R^{-1}$ and 	
	$c' \le \left( \binom{d}{i} - \binom{i_2 - i_1}{i_2 - l} \right) 
              R^{-1}$ $<$ $2^d R^{-1}$.
\item
{\small
\[
 	\variation{i_1}{l}
	\ge
	\log |\zeta_{i_2}| - 2^{-N} \log \left( \binom{i_2-i_1}{i_2 -l} + c'
	\right) + 2^{-N} \log |1+c| 
\]}
	where $|c| \le \left( \binom{d}{i_1} - 1 \right) R^{-1}
	\le 2^d R^{-1}$ and 	
	$c' \le \left( \binom{d}{i} - \binom{i_2 - i_1}{l - i_1} \right) 
              R^{-1}$ $<$ $2^d R^{-1}$. 
\end{enumerate} 
\end{lemma}

\begin{proof}[Proof of Lemma~\ref{cvx3}:]
	By using Lemma~\ref{cvx1} with $i=i_2$, we obtain:
{\small
\begin{equation} \label {ei2}
	-2^{-N} \log | \sigma_{d-i_2} (Z) |
	=
	-2^{-N} \left( \log |Z_{i_2+1}| + \dots \log |Z_d| \right)
	-2^{-N} \log(|1+c''|)
\end{equation}}
	
	Using the same lemma with $i=i_1$, we get:
{\small
\begin{equation} \label {ei1}
	-2^{-N} \log | \sigma_{d-i_1} (Z) |
	=
	-2^{-N} \left( \log |Z_{i_1+1}| + \dots \log |Z_d| \right)
	-2^{-N} \log(|1+c''|)
\end{equation}}
	Subtracting the two previous expressions and dividing by $i_2-i_1$
	we get:
\begin{eqnarray*}
	\variation{i_1}{i_2}
	&=& 2^{-N} \log | Z_{i_2} | + 2^{-N+1} \log(1+c'')\\
	&=& \log | \zeta_{i_2} | + 2^{-N+1} \log(1+c'')
\end{eqnarray*}
	This shows the first part of the Lemma.
\medskip
\par
	By using Lemma~\ref{cvx2}, we can also bound:
\begin{equation}\label{el}
\begin{array}{rcl}
	-2^{-N} \log | \sigma_{d-l} (Z) |
	&\ge&
	-2^{-N} (i_2 - l) \log |Z_{i_2}| \\
        & &
	-2^{-N} \left( \log |Z_{i_2+1}| + \dots \log |Z_d| \right) \\
	& & -2^{-N} \log(|\binom{i_2-i_1}{i_2-l}+c'|)
\end{array}
\end{equation}
	where $c'$ is as in Lemma~\ref{cvx2}.
\medskip
\par
	We can now estimate equation~(\ref{ei2}) minus equation~(\ref{el}),
	altogether divided by $i_2 - l$:
\[
	\variation{l}{i_2}
	\le
	\log |\zeta_{i_2}| + 2^{-N} \log \left( \binom{i_2-i_1}{i_2 -l} + c'
	\right) + 2^{-N} \log |1+c| 
\]
	We can also estimate equation~(\ref{el}) minus equation~(\ref{ei1}),
        altogether divided by $l - i_1$:
\[
        \variation{i_1}{l}
        \ge
        \log |\zeta_{i_2}| - 2^{-N} \log \left( \binom{i_2-i_1}{i_2 -l} + c'
        \right) + 2^{-N} \log |1+c|
\]
\end{proof}

\subsection{A Decision Criterion}
	Lemma~\ref{cvx3} can be used do decide if a point in the 
	convex hull of $g$ is converging to a sharp corner of the limiting
	convex hull or not.

\begin{lemma}\label{cvx4} Assume that
	\begin{enumerate}
	\item [a.] $m \ge \max (i_2 - i_1)$ when $i_1$ and $i_2$ are successive
	 	elements of $I$.
	\item [b.] $2^{-N+1} \log (2^m + 2^d R^{-1}) 
               - 2^{-N+2} \log(1 - 2^d R^{-1}) <
		E < \frac{\log \rho}{2}$
	\item [c.] $i < j < k$
	\end{enumerate}
	Then,
	\begin{enumerate}
	\item If $i$ and $j$ are successive elements of $I$ and there is
	      no other element of $I$ between $j$ and $k$, then
	      $\variation{i}{j} < \variation{j}{k} - E$
	\item If $i$ and $k$ are successive elements of $I$ then
	      $\variation{i}{j} > \variation{j}{k} - E$
	\end{enumerate}
\end{lemma}

\begin{proof}[Proof of Lemma~\ref{cvx4}:]

	Part 1: Assume that $i$, $j$ are successive elements of $I$.
	Then, part 1 of Lemma~\ref{cvx3} implies:
\[
	\variation{i}{j} \le
	\log |\zeta_j| 
	- 2^{-N+1} \log (1 - 2^d R^{-1})
\]
\par
	For the evaluation of $\variation{j}{k}$, we have to
	distinguish two cases: If $k \in I$, then
\[
	\variation{j}{k} \ge
	\log |\zeta_k| 
	+ 2^{-N+1} \log (1 - 2^d R^{-1})
\]
\par
	If $k \not \in I$, let $m$ be such that $j$ and $m$ are 
	successive elements of $I$. Recall that $j<k<m$ by hypothesis.
	Using part~3 of Lemma~\ref{cvx3}, we get:
\[
	\variation{j}{k} \ge
	\log |\zeta_m| 
        - 2^{-N} \log ( 2^m + 2^d R^{-1})
	+ 2^{-N} \log (1 - 2^d R^{-1})
\]
\par
	In any case,
\begin{eqnarray*}
	\variation{i}{j} &\le&
	\variation{j}{k}
	+ \log|\zeta_j| - \log |\zeta_k| \\
	& & + 2^{-N} \log ( 2^m + 2^d R^{-1})
        - 2^{-N+2} \log (1 - 2^d R^{-1})
\end{eqnarray*}
	We use the hypothesis $E < \frac{\log \rho}{2}$ to
	deduce that $\log|\zeta_j| - \log |\zeta_k| + E $ $< -E$,
	and:	
\[
        \variation{i}{j} <
        \variation{j}{k}
	- E
\]
\medskip
\par
	Part 2: Using Lemma~\ref{cvx3}, we have:
\[
\variation{i}{j} \ge
	\log |\zeta_k| 
        - 2^{-N} \log ( 2^m + 2^d R^{-1})
	+ 2^{-N} \log (1 - 2^d R^{-1})
\]
\[
\variation{j}{k} \le
	\log |\zeta_k| 
        + 2^{-N} \log ( 2^m + 2^d R^{-1})
	- 2^{-N} \log (1 - 2^d R^{-1})
\]
\par
	Subtracting, we obtain:
\begin{eqnarray*}
	\variation{i}{j} &\ge&
	\variation{j}{k}
        - 2^{-N+1} \log ( 2^m + 2^d R^{-1})\\
	& & + 2^{-N+1} \log (1 - 2^d R^{-1}) \\
     	&>& \variation{j}{k} - E
\end{eqnarray*} 
\end{proof}

\subsection{Proof of Proposition ~\ref{prop:cvxhull}}
\begin{proof}[Proof of Proposition~\ref{prop:cvxhull}:]

Let $N > 3+\log_2 \frac{ d \log 2 } {\log \rho}$. It is easy
to check that $R > 2^{8d}$, hence $2^d R^{-1} < 2^{-7d}$.
So we can bound: 
\[
	2^{-N+1} \log (2^m + 2^d R^{-1}) 
               - 2^{-N+2} \log(1 - 2^d R^{-1}) 
	<
\]
\[
	<  
	\frac{2(m+1)\log \rho \log 2}{8d \log 2}
	+
	\frac{4 \log \rho}{8 d \log 2} 
	\frac{1}{2^8 - 1} 
	< 
	\frac{\log \rho}{2}
\]

	Therefore, we are in the conditions 1 and 2 of Lemma~\ref{cvx4},
	with $m=d$. Correctness of the algorithm~\ref{alg:cvxhull} can
	be proved now by induction. 
\begin{induction}
	At step $i$, the list $\Lambda$ contains $\Lambda_0, 
	\dots , \Lambda_s$, $\Lambda_{s+1}, \Lambda_j$ where
	$\Lambda_0, \dots \Lambda_s$ are all successive elements
	of $I$ and $\Lambda_{s+1}$, $\dots$, $\Lambda_j$ are not in
	$I$. (Possibly, we can have $s=j$). 
\end{induction}

	The induction hypothesis is true at step 1, with $j=0$, 
        and $0 \in I$. At each step, there are two possibilities:
\par
	Case 1: $i \not \in I$. In that case a few of the
	$\Lambda_{s+1}, \dots \Lambda_j$ may be discarded; but
	part 1 of Lemma~\ref{cvx4} prevents the algorithm from
	discarding elements of $I$.
\par
	Case 2: $i \in I$. In that case, part 2 of Lemma~\ref{cvx4}
	guarantees that all the $\Lambda_{s+1}, \dots \Lambda_j$
	will be discarded. 
\par
	Hence, the induction hypothesis is true at step $i+1$.
	At step $d$, the last point $d$ is added to $\Lambda$.
	Since $d \in I$, $\Lambda = I$.
\medskip
\par
	A note on the running time: although the usual complexity
	of a convex hull algorithm is $O(d \log d)$ for $d$ points
	in the plane, the complexity is smaller when those points
	are `ordered' like ours: $(i, r(i))$. (Compare with Theorem~4.12 in
        ~\cite{PS}).
        Algorithm
	~\ref{alg:cvxhull} has a running time of $O(d)$ operations
	(including a fixed number of transcendental operations). 
	Indeed, each point is added to the list $\Lambda$ precisely
	one time. It can be discarded only once, so the interior
	`while' loop is executed at most $d-1$ times in one
	execution of the algorithm.	
\end{proof}

\section{Tangent Graeffe Iteration}
\label{tangra}
\subsection{Perturbation Methods, Infinitesimals, 1-Jets of Polynomials}

    Graeffe iteration provides the absolute values
    of each root in the case such  roots are 
    all of different moduli. 
Recovering the actual value of each root, and recovering
   pairs of conjugate roots or multiple roots require further work.
\par
   Many algorithms have been proposed to recover the actual roots, such as 
   reverse Graeffe iteration, splitting algorithms. See~\cite{PAN} and
   references therein.  
\par
   A possibility of theoretical interest would be to consider a perturbation of $f$;
   assume first that $f$ is a polynomial with roots $\zeta_1$, \dots , $\zeta_d$ 
   such that $|\zeta_1| < |\zeta_2| < \dots < |\zeta_d|$. 
Then, consider also the
   iterates of 
\[
   f(x+\epsilon) 
\] 
\par
   Graeffe iteration of
   $f(x)$ will provide $|\zeta_1|$, \dots, $|\zeta_d|$, while Graeffe iteration of
   $f(x+\epsilon)$ will provide $|\zeta_1 - \epsilon|$, \dots, $|\zeta_d - \epsilon|$.
   Therefore, we will be able to compute:
\[
   |\zeta_i|^ 2 - |\zeta_i - \epsilon|^2 = - 2 \epsilon \Re(\zeta_i) + \epsilon^2 
\] 
   thus recovering $\zeta_i$.
\medskip
\par
   As mentioned before, this is a possibility of theoretical interest only. The
   perturbation method above would lose half of the working precision in any
   reasonable implementation. Therefore, we will prefer to compute the derivative
   of $|\zeta-\epsilon|^2$ with respect to $\epsilon$. 
\par
   The value $\epsilon$ will be treated as an infinitesimal; therefore, instead of
   storing in memory a certain value $z+ \epsilon \dot z$, we will store 
   $z$ and $\dot z$ separately. When computing some differentiable 
   function $G(z+\epsilon \dot z)$,
   we will obtain a result $G(x) + \epsilon DG(z)\dot z$. So we will compute $G(z)$ and
   $DG(z)\dot z$, but we will never need to assign 
   an actual value to $\epsilon$.
\medskip
\par
   A quantity of the form $z + \epsilon \dot z$ is called a 1-jet. It can also be 
   interpreted as an element of the tangent bundle of the manifold where $z$ is
   supposed to live.
\par
   In this paper, we are mainly concerned with 1-jets of polynomials. We will represent
   degree $d$ polynomials as points in $\mathbb R^{d+1}$ or $\mathbb C ^{d+1}$. Therefore,
   a 1-Jet of polynomials can be represented as a point of $\mathbb R^{2d+2}$ or 
   $\mathbb C ^{2d+2}$, since we are working with a linear space.  
\newtext{

The dot notation (such as in $\dot z$) will be reserved in this paper
to the `tangent' coordinate of a 1-jet $z + \epsilon \dot z$. We reserve the
notation $f'$ to the derivative $\frac{\partial}{\partial x} f$
of a univariate function $f = f(x)$,
and the notation $DF$ to the derivative of a multivariate function $F$.
}{Clarification of the notations for derivatives}
   We need the following construction from Calculus on Manifolds
   ~\cite{ARNOLD,LANG}: Let $G$ be 
   a differentiable function from manifold $X$ into manifold $Y$. Its tangent map
   can be written, in our 1-Jet notation, as:
\[
   \begin{array}{crcl}TG: & T X & \rightarrow & T Y \\ & f + \epsilon \dot f 
 & \mapsto & G(f) + \epsilon DG_f \dot f \mbox{ ,} 

   \end{array} 
\]
   where as usual $TX$ and $TY$ denote the tangent bundle of $X$ and $Y$
   respectively.
\medskip
\par
	The iteration of the $1$-jet of polynomials $f + \epsilon f'$ can
	be used to recover the actual value of the roots of $f$. For
        instance:

\begin{example}
        Let $f$ be a real \condition polynomial, not vanishing at $0$.
	Consider the 1-jet $f(x+\epsilon) = 
        f(x) + \epsilon f'(x)$; 
	its solutions are $\zeta_j - \epsilon$, where $\zeta_j$ are the
	roots of $f$. Let $g + \epsilon \dot g = TG^N (f + \epsilon f')$.
	Lemma~\ref{multiple-pair} below will imply,
        in the particular case $\zeta_j$ is a real isolated root,
        that: 
\[
	\zeta_j = 
        \lim_{N \rightarrow \infty}
	2^{-N}
        \left( \frac{ |g_{j}| }{|g _{j-1}|} \right)
        ^{\frac{1}{2^{N-1}}}
        \left(
        \frac{ \dot g_{j} }{ g_{j} } - \frac{ \dot g_{j-1} }{ g_{j-1} }
        \right)
\]
	\par
	In case $\zeta_j$ and $\zeta_{j+1} = \bar \zeta_j$ are an isolated
	pair of conjugate roots, the limit will be:
\[
	\Re \zeta_j =
        \lim_{N \rightarrow \infty}
	2^{-N-1}
        \left( \frac{ |g_{j+1}| }{|g _{j-1}|} \right)
        ^{\frac{1}{2^{N}}}
        \left(
        \frac{ \dot g_{j+1} }{ g_{j+1} } - \frac{ \dot g_{j-1} }{ g_{j-1} }
        \right)
\]
\end{example}

\par
   In the following section, we compute the tangent map of the Graeffe operator
   in usual and renormalized coordinates. 

\subsection{The Iteration}

	Let $f + \epsilon \dot f$ be a 1-Jet of polynomials. Then its Tangent
     	Graeffe Iterate is:
{\small
\[
	TG ( f(x) + \epsilon \dot f(x) ) = (-1)^d 
        \left( f(\sqrt{x}) + \epsilon \dot f(\sqrt{x}) \right)
        \left( f(-\sqrt{x}) + \epsilon \dot f(-\sqrt{x}) \right)
\]}
\par
	This can be rewritten as:
{\small
\[
	TG ( f(x) + \epsilon \dot f(x) ) = G(f) 
                       + (-1)^d \epsilon \left( f(\sqrt{x}) \dot f(-\sqrt{x}) 
                                              + f(-\sqrt{x}) \dot f(\sqrt{x}) \right) 
\]}
	Precise formulae for computing $g + \epsilon \dot g = TG (f + \epsilon \dot f)$ are:
\[
\left\{
	\begin{array}{ccl}
        g_i      &=& (-1)^{d+i} f_i ^ 2 + 2 \sum_{j \ge 1} (-1)^{d+i+j} f_{i+j} f_{i-j} \\
          \ \\
        \dot g_i &=& 2 \sum_j (-1)^{d+i+j} f_{i-j} \dot f_{i+j} \\ 
        \end{array}
\right.
\]
	For an efficient root-finding algorithm 
	the equations above need to be renormalized.
	At each step, this is done by replacing products and sums 
	by their renormalized counterparts.
        An adjustment is necessary to pass from one renormalization level to
        another (division of the coordinates $r$ by 2). Those adjustments are
	summarized in Algorithm~\ref{alg:tangra}
\begin{algorithm}[!ht]
\caption{TangentGraeffe ($N$, $d$, $r$, $\alpha$, $\hat r$, $\hat \alpha$) }
\label{alg:tangra}
\begin{algorithmic}
\STATE
\COMMENT {~ $N$ (Renormalization level) and $d$ (degree) are integers; $r$ and $\hat r$ should
           be real arrays, and $\alpha$ and $\hat \alpha$ should be array of modulus one complex
           numbers. This routine computes, in renormalized coordinates, the Tangent Graeffe
           Iterate of the 1-jet: $\sum_{i} \left( \alpha_i e^{-2^N r_i} + 
           \epsilon \hat \alpha_i e^{-2^N \hat r_i} \right) x^i$. The coordinates of
           the result are given in renormalization level $N+1$.} 

\STATE $p \leftarrow 2^{N+1}$
\FOR {$i \leftarrow 0$ to $d$}
     \STATE $(s_i, \beta_i)  \leftarrow ( r_i, (-1)^i \alpha_i ^2) $
     \STATE $(\hat s_i, \hat \beta_i )    \leftarrow 
            \left( (r_i + \hat r_i)/2 -  \log(2)/p, 
		(-1)^i \alpha_i \hat \alpha_i \right) $

     \FOR {$j \leftarrow 1$ to $\min(d-i,i)$}

     \STATE $( s_i, \beta_i) \leftarrow 
        \text{RenSum} ( s_i, \beta_i,
                        (r_{i+j} + r_{i-j}) / 2 + \log(2)/p , 
                        (-1)^{i+j} \alpha_{i+j} \alpha_{i-j},p )$
     \STATE $( \hat s_i, \hat \beta_i) \leftarrow 
        \text{RenSum} ( \hat s_i, \hat \beta_i,
                        (r_{i+j} + \hat r_{i-j}) / 2 + \log(2)/p , 
                        (-1)^{i+j} \alpha_{i+j} \hat \alpha_{i-j},p )$
     \STATE $( \hat s_i, \hat \beta_i) \leftarrow 
        \text{RenSum} ( \hat s_i, \hat \beta_i,
                        (r_{i-j} + \hat r_{i+j}) / 2 + \log(2)/p , 
                        (-1)^{i-j} \alpha_{i-j} \hat \alpha_{i+j},p )$
     \ENDFOR
\ENDFOR

\STATE \textbf{return} $(s, \beta, \hat s, \hat \beta)$

\end{algorithmic}

\end{algorithm}

\subsection{Convergence Results}

	It is now time to show convergence of the (Renormalized) Tangent Graeffe Operator.
	Assume one is given a \condition polynomial $f$. Its roots will be ordered as
	$|\zeta_1| \le |\zeta_2| \le \cdots |\zeta_d|$.  
	One can use the (Renormalized) Newton diagram to collect together the roots
        with same moduli. Those will represent single roots, multiple roots, or
        (in the case of real polynomials) pairs of conjugate roots or pairs of
	multiple conjugate roots.
\medskip
\par

\begin{lemma}[Complex case]\label{multiple-root}
	Let $f$ be a complex circle-free polynomial with roots
	$\zeta_1 \ne 0$, \dots, $\zeta_d$ ordered as in 
        Theorem ~\ref{main}. Let
\[
	\rho \bydef \min_{ |\zeta_{i+1}| > |\zeta_i| } 
                        \frac{|\zeta_{i+1}|}{|\zeta_i|}
\]
	and let $g + \epsilon \dot g = (TG)^N (f + \epsilon f')$. 
	Suppose that 
	$j$ and $j+d'$ are successive elements of 
        $I= \{ i: |\zeta_i| < |\zeta_{i+1}| \} \cup \{0 ; d\}$.
	Then,
\[
	\lim_{N \rightarrow \infty}
	- \frac { 2^{-N} }{d'}
        \overline{
	\left(
	\frac{ \dot g_{j+d'} }{ g_{j+d'} } 
             - \frac{ \dot g_{j} }{ g_{j} }
	\right)}
	=
	\frac{\zeta_{j+d'}}{|\zeta_{j+d'}|^2}             \mbox{ .} 
\]
\par
	Furthermore, the error is bounded by:
\[
2^{d+3} \frac{d}{d'} \frac{|\zeta_d|}{|\zeta_1|} \rho^{-2^N} |\zeta_{j+d'}|^{-1}
\]
\end{lemma}
\medskip
\par
	Lemma~\ref{multiple-root} will be proved in Subsection~\ref{proof-of-lemmata}.
\medskip
\par
	Also, real polynomials have usually pairs of conjugate roots;
	they may have pairs with multiplicity. In that case, we can
	show that:

\begin{lemma}[Real case]\label{multiple-pair}
	Let $f$ be a real circle-free polynomial with roots
	$\zeta_1 \ne 0$, \dots, $\zeta_d$ ordered as in 
        Theorem ~\ref{main}. Let
\[
	\rho \bydef \min_{ |\zeta_{i+1}| > |\zeta_i| } 
                        \frac{|\zeta_{i+1}|}{|\zeta_i|}
\]
	and let $g + \epsilon \dot g = (TG)^N f + \epsilon f'$. 
	Suppose that 
	$j$ and $j+d'$ are successive elements of 
        $I= \{ i: |\zeta_i| < |\zeta_{i+1}| \} \cup \{0 ; d\}$.
	Then,
\[
	\lim_{N \rightarrow \infty}
	- \frac { 2^{-N} }{d'}
	\left(
	\frac{ \dot g_{j+d'} }{ g_{j+d'} } 
             - \frac{ \dot g_{j} }{ g_{j} }
	\right)
	=
	\frac{\Re \zeta_{j+d'}}{|\zeta_{j+d'}|^2}             \mbox{ .} 
\]
\par
	Furthermore, the error is bounded by:
\[
2^{d+3} \frac{d}{d'} \frac{|\zeta_d|}{|\zeta_1|} \rho^{-2^N} |\zeta_{j+d'}|^{-1}
\]
\end{lemma}

	The proof of Lemma~\ref{multiple-pair} is also postponed to subsection
	~\ref{proof-of-lemmata}.
	Lemmas~\ref{multiple-root} and~\ref{multiple-pair}
	can be used to recover the roots of a polynomial from the 
	Tangent Graeffe iterates of its $1$-jet:

\begin{algorithm}[!ht]
\caption {RealRecover ( $N,d,I, r, \alpha, \hat r, \hat \alpha$ )} 
\label{alg:rrec}
\begin{algorithmic}
\STATE
\COMMENT {~This procedure attempts to recover the roots of the
          degree $d$ real polynomial $\sum_{i} \alpha_i e^{-2^N r_i} 
          x^i$. The list of sharp corners of its Newton Diagram
          is supposed given in $I=(I_0,\cdots,I_{1 + \text{size}(I)}) $. 
          See Lemma ~\ref{multiple-pair} for a justification~}

\FOR {$k \leftarrow 0$ to Size($I$)}
	\STATE $d' \leftarrow I_{k+1} - I_k$

	\STATE $(b, \beta) \leftarrow \text{RenSum} \left(
                \hat r_{I_{k+1}} - r_{I_{k+1}} ,
                \frac{\hat \alpha_{I_{k+1}}}{\alpha_{I_{k+1}}} ,
                \hat r_{I_k} - r_{I_k},
                -\frac{\hat \alpha_{I_k}}{\alpha_{I_k}}, 
                2^N
                \right)$
	\STATE
	\STATE $m \leftarrow \exp \left(2 \frac{r_{I_{k+1}} - r_{I_k}}{d'} \right)$ 

	\STATE $x \leftarrow -\beta \frac{2^{-N}}{d'} m \exp -2^N b$ 
        \IF {$I_{k+1}-I_k$ is even and $m > |x|^2$} 
            \STATE $y \leftarrow \sqrt{m - |x|^2}$
	    \ELSE  
	    \STATE $x\leftarrow m\frac{x}{|x|}$
	    \STATE $y\leftarrow 0$
	    \ENDIF

        \FOR {$j \leftarrow 0$ to $I_{k+1} - I_k -1$}
             \STATE $\zeta_{I_k + j +1} = x + (-1)^j y$
        \ENDFOR
\ENDFOR
\STATE \textbf{return} $\zeta$
\end{algorithmic}
\end{algorithm}

\begin{algorithm}[!ht]
\caption {ComplexRecover ( $N,d,I, r, \alpha, \hat r, \hat \alpha$ )} 
\label{alg:crec}
\begin{algorithmic}
\STATE
\COMMENT {~This procedure attempts to recover the roots of the
          degree $d$ complex polynomial $\sum_{i} \alpha_i e^{-2^N r_i} 
          x^i$. The list of sharp corners of its Newton Diagram
          is supposed given in $I$. See Lemma ~\ref{multiple-root} 
          for a justification~}

\FOR {$k \leftarrow 0$ to Size($I$)}
	\STATE $d' \leftarrow I_{k+1} - I_k$

	\STATE $(b, \beta) \leftarrow \text{RenSum} \left(
                \hat r_{I_{k+1}} - r_{I_{k+1}} ,
                \frac{\hat \alpha_{I_{k+1}}}{\alpha_{I_{k+1}}} ,
                \hat r_{I_k} - r_{I_k},
                -\frac{\hat \alpha_{I_k}}{\alpha_{I_k}} 
                2^N,
                \right)$

	\STATE
	\STATE $m \leftarrow \exp \left(2 \frac{r_{I_{k+1}} - r_{I_k}}{d'} \right)$ 

	\STATE $x \leftarrow - \bar \beta \frac{2^{-N}}{d'} m \exp -2^N b$ 

        \FOR {$j \leftarrow 0$ to $I_{k+1} - I_k -1$}
             \STATE $\zeta_{I_k + j +1} = x $
        \ENDFOR
\ENDFOR
\STATE \textbf{return} $\zeta$
\end{algorithmic}
\end{algorithm}

\subsection{The Main Algorithm}

	We can now state the algorithm of Theorem~\ref{main}. 
We start with
a fixed, arbitrary value for $\rho(f) = \max_{|\zeta_i| > |\zeta_j|} 
\frac{|\zeta_i|}{|\zeta_j|}$. 
 Proposition~\ref{prop:cvxhull} guarantees that if 
 \[ N > 3 + \log_2 \frac{d \log 2}{\log \rho(f)} \mbox{ ,}  \] 
then after the $N$-th iterate
the convex hull of the Newton Diagram of $f$ is computed correctly.

\begin{algorithm}[!ht]
\caption {Solve ($d,f, \text{isreal}$)}
\label{alg:solve}

\begin{algorithmic}
\STATE
\COMMENT {~It is assumed here that $f$ is a degree $d$, circle-free real or complex 
          polynomial. In the general case, one should first
          find and output the trivial ($0$ and $\infty$) roots of $f$, then deflate $f$.
          After that, one should perform a random real (resp. complex) conformal transform 
          on $f$ so it becomes circle-free~}
\STATE
\FOR {$i \leftarrow 0$ to $d$}
	\IF    {$f_i \ne 0$}
	\STATE $\alpha_i \leftarrow f_i / |f_i|$
	\ELSE
	\STATE $\alpha_i \leftarrow 1$
	\ENDIF
	\STATE $r_i \leftarrow -\log |f_i|$
	\ENDFOR
\FOR {$i \leftarrow 0$ to $d-1$}
        \STATE $f'_i \leftarrow (i+1) f_{i+1}$
	\IF    {$f'_i \ne 0$}
	\STATE $\hat \alpha_i \leftarrow f'_i / |f'_i|$
	\ELSE
	\STATE $\hat \alpha_i \leftarrow 1$
	\ENDIF
	\STATE $\hat r_i \leftarrow -\log |f'_i|$
	\ENDFOR
\STATE $N \leftarrow 0$
\STATE $\rho = 2$
\LOOP
\STATE $r$, $\alpha, \hat r, \hat \alpha\leftarrow$ 
       TangentGraeffe ($N,d,r, \alpha, \hat r, \hat \alpha$)
\STATE $N \leftarrow N + 1$
\STATE $I \leftarrow$ Convex Hull ($N$, $d$, $r$, $\rho$)
\IF {isreal}
      \STATE $\zeta \leftarrow$ RealRecover ( $N,d,I, r, \alpha, \hat r, \hat \alpha$ )
      \ELSE
      \STATE $\zeta \leftarrow$ ComplexRecover ( $N,d,I, r, \alpha, \hat r, \hat \alpha$ )
      \ENDIF
\STATE Output $\zeta_1, \cdots \zeta_d$
\IF    {$N > 3 + \log_2 \frac{d \log 2}{\log \rho}$}
       \STATE
       \COMMENT {~ Proposition~\ref{prop:cvxhull} implies that at this point,
                $I$ is indeed correct for all the polynomials
                with separation ration $\ge \rho$. Therefore, it is
                time to decrease $\rho$.~}
       \STATE $\rho \leftarrow \sqrt{\rho}$
       \ENDIF
\ENDLOOP

\end{algorithmic}
\end{algorithm}

After the $N$-th iteration, convergence is guaranteed by the following
bounds: According to Lemma~\ref{cvx3}, at the execution of algorithm
Complex Recover (resp. Real Recover),
\[
\left|\exp \left( \frac{2}{i_{k+1}-i_k} (g_{i_{k+1}} - g_{i_k}) \right)\right|
= |\zeta_{i_{k+1}}| (1 + \delta_1)
\]
	where $|\delta_1| \le e^{2^{-N+1} \log 1+ 2^d \rho^{-2^N}}$.
\par
	Introducing the error bound of Lemma~\ref{multiple-root},
one gets in the complex case:
\[
- \frac{2^{-N}}{d'} \overline{(a - b)} = 
\frac{\zeta_{i_{k+1}}} {|\zeta_{i_{k+1}}|^2} 
(1 + \delta_2)
\]
with $|\delta_2| < 
2^{d+3} \frac{d}{d'} \frac{|\zeta_d|}{|\zeta_1|} \rho^{-2^N}$,
and where $a$ and $b$ are as in the Algorithms. 
Therefore,
\[
|\zeta_{i_k + j +1} - \bar x|
\le 
\delta |\zeta_{i_k + j +1}| 
\] 
	where $\delta < (1 + \delta_1)^2 (1 + \delta_2) -1$.

The real case is analogous. According to Lemma ~\ref{multiple-pair},   
\[
- \frac{2^{-N}}{d'} {(a - b)} = 
\frac{\Re \zeta_{i_{k+1}}} {|\zeta_{i_{k+1}}|^2} 
(1 + \delta_2)
\]
\par
so that
\[
|\Re \zeta_{i_k + j +1} - x| < \delta |\zeta_{i_k + j +1}|
\] 
and
\[
|\Im \zeta_{i_k + j +1} - \sqrt{1-x^2}| < \delta' |\zeta_{i_k + j +1}|
\]
\par
where $\delta' = \delta + O(\delta^2)$ 
\medskip
\par
  In both cases, $|\delta|$ and eventually $|\delta'|$ are
dominated by $\rho^{-2^N}$. Since 
\[ A \rho^{-2^N} = 2 ^{\log_2 A - 2^N \log_2 \rho} 2^{-2^{N-C}} \mbox{ ,} \]
the bound in Theorem~\ref{main} follows.

\subsection{Proof of Lemmas ~\ref{multiple-root} and ~\ref{multiple-pair}}
\label{proof-of-lemmata}

	Consider the 1-jet of degree $d$ polynomials $f + \epsilon \dot f$, 
	with
	solutions $\zeta_i + \epsilon \dot \zeta_i$. After $N$ steps
	of Tangent Graeffe Iteration, we obtain a 1-jet of polynomials
\[
	g + \epsilon \dot g = \left( TG \right)^N (f + \epsilon \dot f)  \mbox{ .} 
\]
 Since 
 the differential of the transformation $P^{N}: \zeta \longmapsto z=\zeta^{2^{N}}$ gives
\[ 	DP^{N}: \zeta+\epsilon \dot{\zeta} \longmapsto \left( \zeta_j \right)^{2^N}
      + \epsilon 2^N \left( \zeta_j \right)^{2^N}
        \frac {\dot \zeta_j}{\zeta_j}
\mbox{ ,}  \]
it is clear that the roots $Z_{j}+\epsilon \dot{Z}_{j}$ of $g + \epsilon \dot{g}$
       will be 
\[
        \left( \zeta_j \right)^{2^N} 
      + \epsilon 2^N \left( \zeta_j \right)^{2^N}
	\frac {\dot \zeta_j}{\zeta_j} 	 \mbox{ .}
\]
We can now
compute the derivative at $\epsilon=0$ of 
$g_i + \epsilon g_i = \sigma_{d-i}(Z+\epsilon \dot{Z})$.

Let's denote by $Z_{\widehat{j}}$ the vector $(Z_{1},\cdots,\widehat{Z_{j}},\cdots,Z_{d})
\in {\mathbb C}^{d-1} $.

Take $i \in I$ and notice that
\begin{eqnarray*}
 \frac{\partial}{\partial \epsilon} \restricted{\epsilon=0}  \sigma_{d-i}(Z+\epsilon \dot{Z}) & = & 
\sum_{j} \dot{Z}_{j} \sigma_{d-i-1}(Z_{\widehat{j}}) \\
	& = & \sum_{j} \frac{ \dot{Z}_{j} }{Z_{j}} Z_{j} \sigma_{d-i-1}(Z_{\widehat{j}})
\end{eqnarray*}
Thus,
\begin{eqnarray*}
\frac{\partial}{\partial \epsilon} \restricted{\epsilon=0}
\frac{ \sigma_{d-i}(Z+\epsilon \dot{Z})}{\sigma_{d-i}(Z)} & = & 
\sum_{j} \frac{ \dot{Z}_{j} }{Z_{j}} \frac{Z_{j} \sigma_{d-i-1}(Z_{\widehat{j}})}{\sigma_{d-i}(Z)} \\
& = & \left(\sum_{j>i} + \sum_{j\le i}\right) 
\frac{ \dot{Z}_{j} }{Z_{j}} \frac{Z_{j} \sigma_{d-i-1}(Z_{\widehat{j}})}{\sigma_{d-i}(Z)} 
\mbox{ .} 
\end{eqnarray*}
Due to Lemma~\ref{cvx1}
\begin{eqnarray*}
\sum_{j>i}  \frac{ \dot{Z}_{j} }{Z_{j}} \frac{Z_{j} \sigma_{d-i-1}(Z_{\widehat{j}})}{\sigma_{d-i}(Z)} 
& = &
\sum_{j>i} \frac{\dot{Z}_{j}}{Z_{j}} (1 + \eta_{j}) 
\end{eqnarray*}
and
\begin{eqnarray*}
\sum_{j\le i}  \frac{ \dot{Z}_{j} }{Z_{j}} \frac{Z_{j} \sigma_{d-i-1}(Z_{\widehat{j}})}{\sigma_{d-i}(Z)} 
& = &
\sum_{j\le i} \frac{\dot{Z}_{j}}{Z_{j}} \frac{Z_{j}}{Z_{i+1}} 
 \left(1 + \eta_{j} \right) 
\mbox{ ,} 
\end{eqnarray*}
where  
\[ |\eta_{j}| \le 4 \binom{d}{i} R^{-1} \mbox{ .} \]
Although this estimate is by no means sharp, it suffices for our purposes.
\medskip
\par
Using that $i\in I$ and hence $|z_{j}/z_{i+1}|< R^{-1}$ for $i\le j$, we get 
\begin{small}
\begin{eqnarray*}
\Bigg| \frac{\partial}{\partial \epsilon} \restricted{\epsilon=0}  
\frac{ \sigma_{d-i}(z+\epsilon \dot{z})}{\sigma_{d-i}(z)}   - 
\sum_{j>i} \frac{\dot{z}_{j}}{z_{j}} \Bigg|  & \le  & 
\max \bigg|    \frac{\dot{z}_{j}}{z_{j}} \bigg| 
 \frac{4 \binom{d}{i}(d-i) + i \left( 1 + 4 \binom{d}{i} R^{-1} \right)}{R} \\
	& < & \max\bigg|    \frac{\dot{z}_{j}}{z_{j}} \bigg| d  2^{d+2} R^{-1} 
\end{eqnarray*}
\end{small}

Therefore, if we take the logarithmic derivative of the expression  
\[ (g + \epsilon \dot{g}) = (TG)^{N} ( f + \epsilon \dot{f} ) \]
and evaluate at a successive couple of elements $i_{1},i_{2}\in I$
we get
\[ \bigg| 
          \left( \frac{\dot{g}_{i_{1}}}{g_{i_{1}}}  
          - \frac{\dot{g}_{i_{2}}}{g_{i_{2}}} \right)
+ \sum_{i_1 < j  \le i_2} 2^{N} \frac{1}{\zeta_{j}} \bigg| < 
\frac{2^{N+d+3} d}{R} \max_{l}  \bigg|  \frac{1}{\zeta_{l}}  \bigg| \mbox{ .} \]

Now let's assume we are under the hypothesis of Lemma~\ref{multiple-root}.
Then, since $i_1$ and $i_2$ are consecutive, it follows that
\[ \sum_{i_1 < j \le i_2} \zeta_{j}^{-1} = 
\sum_{i_1 < j \le i_2} \frac{\overline{\zeta_{j}}}{| \zeta_{j}|^{2}} = 
(i_2 -i_1) \frac{\overline{\zeta_{i_2}}}{| \zeta_{i_2}|^{2}} \]
and so
\begin{equation} \label{compcase}
\bigg| \frac{2^{-N}}{i_2-i_1} \left(\frac{\dot{g}_{i_{2}}}{g_{i_{2}}}   - 
\frac{\dot{g}_{i_{1}}}{g_{i_{1}}} \right) |\zeta_{i_2}|^{2} 
-\overline{\zeta_{i_{2}}} \bigg| |\zeta_{i_2}|^{-1}
< \frac{ 2^{d+3} d}{ R (i_2 -i_1) }  \max_{r,s} \frac{|\zeta_{r}|}{|\zeta_{s}|}
\mbox{ .} 
\end{equation}

On the other hand, if we are under the hypothesis of Lemma~\ref{multiple-pair},
we get 
\[ \sum_{i_1 < j \le i_2} \zeta_{j}^{-1} = 
(i_2 - i_1) \frac{\Re{ \zeta_{i_2} }}{|\zeta_{i_2}|^2} \mbox{ .} \]
Thus, 
\begin{equation} \label{realcase}
\bigg| \frac{2^{-N}}{i_2-i_1} \left(\frac{\dot{g}_{i_{2}}}{g_{i_{2}}}   - 
\frac{\dot{g}_{i_{1}}}{g_{i_{1}}} \right) |\zeta_{i_2}|^{2}
-\Re{\zeta_{i_{2}}} \bigg| |\zeta_{i_2}|^{-1}
< \frac{ 2^{d+3} d}{ R (i_2 -i_1) }  \max_{r,s} \frac{|\zeta_{r}|}{|\zeta_{s}|}
\mbox{ .} 
\end{equation}

In either case, we have that each  right hand side
 of equations~(\ref{compcase}) and (\ref{realcase}) is bounded by 
\[ 2^{d+3} \frac{d}{R(i_2 - i_1)}  \frac{|\zeta_{d}|}{|\zeta_{1}|} \mbox{ .}\]
Now, using the fact that the separation radius at the $N$-th step
\[ R = \rho^{2^N} \mbox{ ,} \] 
where $\rho$ is the separation radius of the original roots of $f$ before 
applying Graeffe, we get that
\[  \frac{ 2^{d+3} d}{ R (i_2 -i_1) } \frac{|\zeta_{d}|}{|\zeta_{1}|} =
2^{d+3}\frac{d}{i_2 - i_1}  \frac{|\zeta_{d}|}{|\zeta_{1}|} \rho^{-2^N}
\longrightarrow 0 \mbox{ ,  as } N \longrightarrow \infty  \mbox{ .}  \]
\myqed

\newtext{
\subsection{`Deterioration of Condition' and Stability Properties}
\label{detcond}
It is important to understand that we will never have to solve $g(x)=$
$(G^N f)(x)$. Therefore, the actual condition number of $g$ does
not matter at all. In order to determine the roots of $f$, we
will be using the extra information provided by $\dot g$, where
$(g,\dot g) = TG^N (f,\dot f)$. 

A valid source of concern is the propagation of rounding-off error.
In the Tangent Graeffe algorithm, that error would typically 
double at each step
(it actually doubles at each step `in the limit'). 

However, Lemmas~\ref{multiple-root} and ~\ref{multiple-pair} guarantee that the truncation 
error decreases as $\rho^{-2^N}$. Hence, in order to obtain 
a truncation error smaller than a certain $\delta > 0$, we need
$N \approx c + \log_2 \log_2 \delta^{-1}$, where $c=-\log_2 \log_2 \rho$ 
is a constant depending only on the original polynomial $f$. 

In order to reduce the accumulated rounding-off error to the same
order, one would need that
\[
C \macheps 2^N < \delta
\]
where $\macheps$ is the `machine epsilon' and $C$ is a constant
depending on $f$. Thus, we just need
\[
\macheps 
<
\frac{\delta}{C}2^{-N}
\approx
\frac{\delta}{C 2^c \log_2 \delta^{-1}}
\]

What is a reasonable value for $\delta$ ? The strength of
Renormalized Tangent Graeffe Iteration is its capacity to
solve the `global' problem: given a polynomial, approximate
all its roots. Once a suitable approximation of each root
is found, local iterative algorithms (such as Newton Iteration)
cheaply provide better refinements of the roots. Such a two-step
procedure entails a reasonable range of values for $\delta$.
Namely, $\delta$ should be smaller (but not much smaller) than the
radius of quadratic convergence of Newton's iteration.

This radius is of the order of the reciprocal condition number
of the original polynomial.
(See ~\cite{BCSS} Theorem~1 and Remark~1 p. 263 
for a precise statement).

}{Here we answer to the referee's concern on the stability of
this iterative process}
\section{Numerical Results and Final Remarks}
\label{numerics}
\subsection{Numerical Results}

	A polynomial solver based on the algorithms above
	was implemented and tested under IEEE 754 double and
	double-extended arithmetic. The results below are 
	intended to make a case in favor of the stability and
	the practical feasibility of Renormalized Tangent Graeffe
	Iteration.

\begin{figure}
\centerline{\epsfig{file=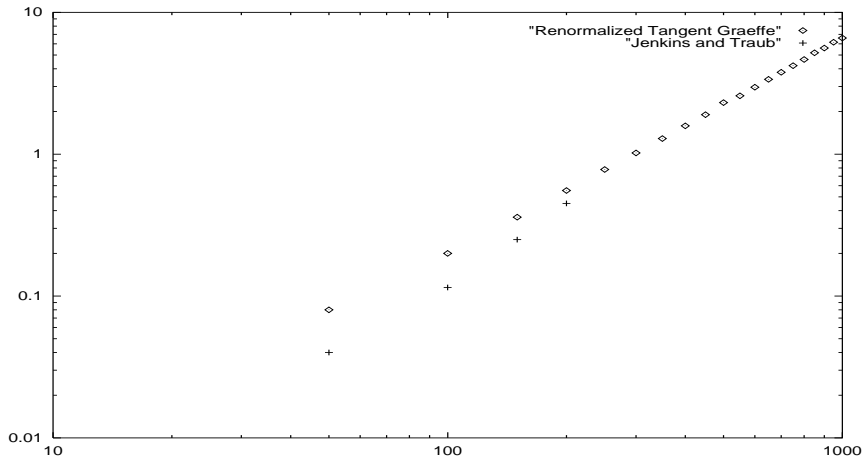, height=6cm, width=12cm}}
\caption[ ]{\label{fig1}Real pseudo-random polynomials}
\end{figure}

\begin{figure}
\centerline{\epsfig{file=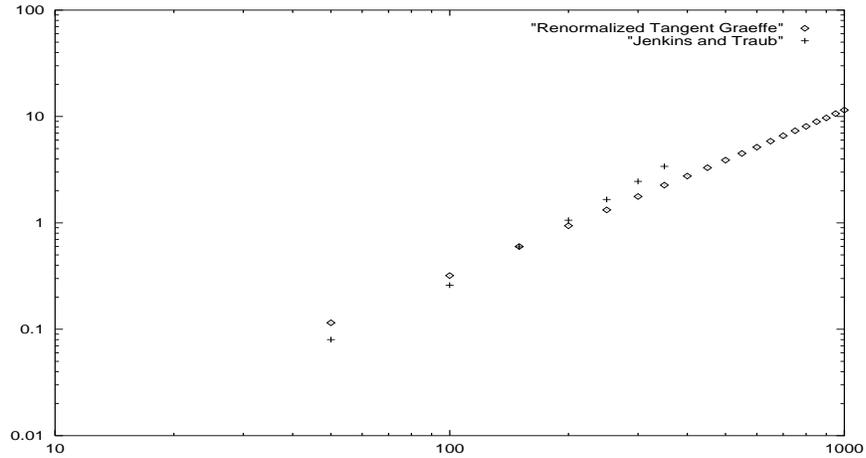, height=6cm, width=12cm}}
\caption[ ]{\label{fig2}Complex pseudo-random polynomials}
\end{figure}

	The first set of tests was designed to measure the performance
	of our algorithm for large degree polynomials. The test polynomials
	are pseudo-random real (Table~\ref{tab1} and Figure~\ref{fig1}) 
        and complex (Table~\ref{tab2} and Figure~\ref{fig2})  
        polynomials, under the
	$U(2)$-invariant probability measure~\cite{KOSTLAN,BEZII,MZ97}. 
        Under this probability measure, random polynomials are
	well-conditioned on the average.

	The results were certified using alpha-theory~\cite{SMALE,BEZI,TCS}. 
        The running time (certification excluded) was 
	compared to the code of Jenkins and \linebreak Traub 
        ~\cite{TOMS419,TOMS493} for the values where this code succeeds.

	Running time is measured in user-time seconds of a Pentium-133 computer
	running Linux and the gcc compiler.
\medskip
\par

\begin{table}
\centerline{\tiny
\begin{tabular}{|r|l|rrrrrrrrrr|}
\hline
\hline
\multicolumn{12}{|c|}{Running time (s)}\\
\hline
\hline
Degree & Algorithm & \multicolumn{10}{c|}{Seed}\\
& & 0 & 1 & 2 & 3 & 4 & 5 & 6 & 7 & 8 & 9 \\
\hline
50 & Graeffe &
0.08&
0.08&
0.08&
0.07&
0.09&
0.08&
0.08&
0.08&
0.08&
0.07\\
 & J-T &
0.03&
0.02&
0.04&
0.02&
0.05&
0.05&
0.07&
0.04&
0.04&
0.03\\
\hline
100 & Graeffe &
0.20&
0.20&
0.20&
0.21&
0.21&
0.21&
0.19&
0.20&
0.19&
0.20\\
 & J-T &
0.12&
0.11&
0.12&
0.09&
0.12&
0.12&
0.14&
0.09&
0.10&
0.10\\
\hline
150 & Graeffe &
0.36&
0.36&
0.36&
0.36&
0.37&
0.35&
0.36&
0.37&
0.35&
0.36\\
 & J-T &
0.26&
0.26&
0.28&
0.24&
0.23&
0.24&
0.27&
0.20&
0.26&
0.24\\
\hline
200 & Graeffe &
0.56&
0.54&
0.56&
0.56&
0.55&
0.55&
0.56&
0.55&
0.56&
0.55\\
 & J-T &
0.53&
0.42&
0.51&
0.35&
0.43&
0.38&
0.54&
0.40&
0.47&
0.36\\
\hline
250 & Graeffe &
0.77&
0.78&
0.77&
0.78&
0.78&
0.78&
0.78&
0.78&
0.77&
0.76\\
\hline
300 & Graeffe &
1.02&
1.02&
1.01&
1.02&
1.02&
1.00&
1.02&
1.01&
1.01&
1.03\\
\hline
350 & Graeffe &
1.29&
1.29&
1.29&
1.29&
1.30&
1.29&
1.29&
1.29&
1.28&
1.29\\
\hline
400 & Graeffe &
1.59&
1.59&
1.57&
1.58&
1.58&
1.60&
1.57&
1.59&
1.59&
1.58\\
\hline
450 & Graeffe &
1.90&
1.90&
1.89&
1.91&
1.89&
1.89&
1.90&
1.90&
1.91&
1.91\\
\hline
500 & Graeffe &
2.32&
2.27&
2.31&
2.30&
2.31&
2.31&
2.30&
2.29&
2.32&
2.34\\
\hline
550 & Graeffe &
2.57&
2.58&
2.58&
2.58&
2.59&
2.59&
2.58&
2.59&
2.57&
2.59\\
\hline
600 & Graeffe &
2.96&
2.98&
2.97&
2.96&
2.95&
2.97&
2.96&
2.97&
2.98&
2.96\\
\hline
650 & Graeffe &
3.38&
3.39&
3.37&
3.36&
3.37&
3.34&
3.38&
3.37&
3.37&
3.37\\
\hline
700 & Graeffe &
3.78&
3.78&
3.78&
3.80&
3.79&
3.77&
3.78&
3.80&
3.78&
3.77\\
\hline
750 & Graeffe &
4.21&
4.22&
4.21&
4.21&
4.19&
4.20&
4.19&
4.22&
4.20&
4.21\\
\hline
800 & Graeffe &
4.66&
4.66&
4.65&
4.65&
4.64&
4.64&
4.63&
4.65&
4.65&
4.65\\
\hline
850 & Graeffe &
5.18&
5.20&
5.18&
5.17&
5.23&
5.18&
5.19&
5.20&
5.19&
5.17\\
\hline
900 & Graeffe &
5.61&
5.61&
5.58&
5.59&
5.60&
5.57&
5.56&
5.60&
5.60&
5.60\\
\hline
950 & Graeffe &
6.16&
6.16&
6.15&
6.16&
6.13&
6.14&
6.16&
6.18&
6.14&
6.18\\
\hline
1000 & Graeffe &
6.60&
6.58&
6.58&
6.58&
6.59&
6.61&
6.59&
6.60&
6.58&
6.59\\
\hline
\end{tabular}
}
\caption{Real pseudo-random polynomials \label{tab1}}
\end{table}

\begin{table}
\centerline{\tiny
\begin{tabular}{|r|l|rrrrrrrrrr|}
\hline
\hline
\multicolumn{12}{|c|}{Running time (s)}\\
\hline
\hline
Degree & Algorithm & \multicolumn{10}{c|}{Seed}\\
& & 0 & 1 & 2 & 3 & 4 & 5 & 6 & 7 & 8 & 9 \\
\hline
50 & Graeffe &
0.13&
0.12&
0.12&
0.13&
0.12&
0.10&
0.11&
0.10&
0.11&
0.11\\
& J-T &
0.07&
0.08&
0.09&
0.09&
0.07&
0.08&
0.07&
0.08&
0.07&
0.09\\
\hline
100 & Graeffe &
0.34&
0.32&
0.32&
0.32&
0.33&
0.33&
0.30&
0.33&
0.31&
0.32\\
& J-T &
0.25&
0.28&
0.26&
0.25&
0.28&
0.26&
0.26&
0.27&
0.26&
0.25\\
\hline
150 & Graeffe &
0.60&
0.60&
0.61&
0.60&
0.61&
0.60&
0.60&
0.61&
0.60&
0.60\\
& J-T &
0.55&
0.60&
0.58&
0.66&
0.61&
0.60&
0.60&
0.59&
0.59&
0.59\\
\hline
200 & Graeffe &
0.95&
0.95&
0.97&
0.94&
0.96&
0.94&
0.92&
0.94&
0.94&
0.93\\
& J-T &
1.13&
1.17&
1.07&
1.04&
1.09&
1.06&
1.05&
1.06&
1.03&
1.05\\
\hline
250 & Graeffe &
1.33&
1.34&
1.30&
1.32&
1.35&
1.32&
1.32&
1.32&
1.35&
1.33\\
&J-T&
1.63&
1.68&
1.64&
1.68&
1.66&
1.70&
1.66&
1.63&
1.63&
1.66\\
\hline
300 & Graeffe &
1.77&
1.76&
1.77&
1.78&
1.77&
1.76&
1.77&
1.75&
1.79&
1.76\\
&J-T&
2.34&
2.49&
2.43&
2.44&
2.50&
2.46&
2.48&
2.45&
2.36&
2.51\\
\hline
350 & Graeffe &
2.26&
2.30&
2.26&
2.25&
2.28&
2.34&
2.08&
2.36&
2.27&
2.26\\
&J-T&
3.30&
3.41&
3.39&
3.34&
3.36&
3.72&
3.51&
3.45&
3.35&
3.52\\
\hline
400 & Graeffe &
2.74&
2.75&
2.77&
2.76&
2.77&
2.76&
2.75&
2.77&
2.74&
2.74\\
\hline
450 & Graeffe &
3.27&
3.28&
3.28&
3.30&
3.38&
3.30&
3.30&
3.28&
3.34&
3.30\\
\hline
500 & Graeffe &
3.92&
3.87&
3.89&
3.89&
3.88&
3.89&
3.91&
3.87&
3.90&
3.91\\
\hline
550 & Graeffe &
4.49&
4.48&
4.51&
4.52&
4.50&
4.52&
4.51&
4.49&
4.50&
4.48\\
\hline
600 & Graeffe &
5.16&
5.14&
5.17&
5.18&
5.15&
5.14&
5.15&
5.15&
5.18&
5.13\\
\hline
650 & Graeffe &
5.89&
5.88&
5.87&
5.83&
5.83&
5.86&
5.89&
5.84&
5.89&
5.86\\
\hline
700 & Graeffe &
6.62&
6.59&
6.59&
6.58&
6.59&
6.61&
6.59&
6.59&
6.58&
6.62\\
\hline
750 & Graeffe &
7.37&
7.42&
7.40&
7.32&
7.48&
7.27&
7.30&
7.30&
7.43&
7.31\\
\hline
800 & Graeffe &
8.10&
8.06&
8.07&
8.03&
8.10&
8.07&
8.10&
8.06&
8.11&
8.09\\
\hline
850 & Graeffe &
8.96&
8.94&
8.92&
8.95&
8.97&
8.94&
8.92&
8.92&
8.95&
8.96\\
\hline
900 & Graeffe &
9.72&
9.70&
9.70&
9.70&
9.72&
9.72&
9.71&
9.70&
9.67&
9.71\\
\hline
950 & Graeffe &
10.63&
10.67&
10.65&
10.64&
10.67&
10.65&
10.67&
10.61&
10.65&
10.60\\
\hline
1000 & Graeffe &
11.49&
11.47&
11.54&
11.46&
11.50&
11.49&
11.54&
11.50&
11.51&
11.50\\
\hline
\end{tabular}
}
\caption{Complex pseudo-random polynomials
\label{tab2}}
\end{table}

\medskip
\par

	In the second set of experiments, we tried to check the 
	behavior of Renormalized Tangent Graeffe Iteration in the 
	presence of very badly conditioned polynomials. The test
	polynomials are Wilkinson' s `perfidious' polynomials~\cite{PERF} 
\[
	p_d(x) = (x-1)(x-2) \cdots (x-d)
\]

	and Chebyshev polynomials (Table~\ref{tab3}).
\[
	T_d = \prod_{0 \le m < d}
	\left(x - \cos ( \frac{\pi}{2 d} + \frac{\pi}{d}m ) \right)
\]

	The error of the solutions of the perfidious polynomials
	is measured as $\max |\zeta - \text{round \ } \zeta |$, $\zeta$ a
	`solution' found by the program. Similarly, the error in
	Chebyshev polynomials is measured as $\max |m - \text{ round \ } m|$ where
	$m = \frac{d \arccos \zeta - \frac{\pi}{2}}{\pi}$ and
	$\zeta$ a `solution' found by the program. Again, those 
	results are compared to the ones provided by the software
	by Jenkins and Traub. 

\begin{table}
\centerline{\tiny
\begin{tabular}{|r|l|r|}
\hline
\hline
\multicolumn{3}{|c|}{Perfidious polynomials}\\
\hline
\hline
Degree & Algorithm & \multicolumn{1}{c|}{Error}\\
\hline
10 & Graeffe & $5.123013 \times 10^{-12}$ \\
   & J-T     & $4.859935 \times 10^{-11}$ \\
\hline
15 & Graeffe & $3.968295 \times 10^{-08}$ \\
   & J-T     & $5.508868 \times 10^{-09}$ \\
\hline
20 & Graeffe & $1.780775 \times 10^{-03}$ \\
   & J-T     & $1.275754 \times 10^{-04}$ \\
\hline
\end{tabular}
\hspace{1cm}
\begin{tabular}{|r|l|r|}
\hline
\hline
\multicolumn{3}{|c|}{Chebyshev polynomials}\\
\hline
\hline
Degree & Algorithm & \multicolumn{1}{c|}{Error}\\
\hline
10 & Graeffe & $8.790711\times 10^{-16}$\\
   & J-T     & $8.790711\times 10^{-16}$\\
\hline
15 & Graeffe & $2.169163\times 10^{-15}$\\
   & J-T     & $2.169163\times 10^{-15}$\\
\hline
20 & Graeffe & $1.903848\times 10^{-14}$\\
   & J-T     & $1.278977\times 10^{-13}$\\
\hline
25 & Graeffe & $1.266375\times 10^{-11}$\\
   & J-T     & $1.663025\times 10^{-11}$\\
\hline
30 & Graeffe & $5.511325\times 10^{-11}$\\
   & J-T     & $3.301608\times 10^{-10}$\\
\hline
35 & Graeffe & $5.708941\times 10^{-09}$\\
   & J-T     & $3.840000\times 10^{-08}$\\
\hline
\end{tabular}
}
\caption{Perfidious and Chebyshev polynomials \label{tab3}}
\end{table}

\subsection{Further Practical Remarks}

\begin{itemize} 
\item Graeffe process (and hence our algorithm) is known to be 
       parallelizable~ \cite{JS,RJ}.
\item The algorithm presented here needs only $O(d)$ memory
      storage; therefore, all intermediate computations for
      a reasonable degree $d$ fit into the `cache' memory of
      modern computers.  
\item Use of higher precision may be required to handle very
      badly conditioned polynomials; as a matter of fact, 
      polynomials of that sort are only meaningful if their coefficients
      are known with sufficiently high accuracy. For example, when they
	are  obtained by symbolic manipulation.  

\end{itemize}

\section*{Acknowledgements}
JPZ's work was supported in part by CNPq, through grant MA \linebreak 521.329/94-9, 
and by PRONEX grant 76.97.1008-00. GM's work was supported by CNPq through 
grant MA 520.423/96-8 and FAPERJ E-26/170.027/98.



\begin{thebibliography}{10}

\bibitem{ARNOLD}
{\sc V.~Arnold}, {\em M{\'e}thodes Math{\'e}matiques de la M{\'e}canique
  Classique}, Mir, Moscow, 1976.

\bibitem{BP}
{\sc D.~Bini and V.~Y. Pan}, {\em Graeffe's, {C}hebyshev-like, and {C}ardinal's
  processes for splitting a polynomial into factors}, J. Complexity, 12 (1996),
  pp.~492--511.
\newblock Special issue for the Foundations of Computational Mathematics
  Conference (Rio de Janeiro, 1997).

\bibitem{BCSS}
{\sc L.~Blum, F.~Cucker, M.~Shub, and S.~ Smale},
\newblock {\em Complexity and Real Computation}.
\newblock Springer, 1998.

\bibitem{BS}
{\sc Mr S.~Brodetsky and G.~Smeal}, {\em On {G}raeffe's method for
complex roots of algebraic equations}, Proceedings of the Cambridge
Philosophical Society XXII Part II (1924), pp.~83--87

\bibitem{DEDIEU}
{\sc J.-P. Dedieu}, {\em \`{A} propos de la m\'ethode de {D}andelin-{G}raeffe},
  C. R. Acad. Sci. Paris S\'er. I Math., 309 (1989), pp.~1019--1022.

\bibitem{DGY}
{\sc J.~P. Dedieu, X.~Gourdon, and J.~C. Yakoubsohn}, {\em Computing the
  distance from a point to an algebraic hypersurface}, in The mathematics of
  numerical analysis (Park City, UT, 1995), vol.~32 of Lectures in Appl. Math.,
  Amer. Math. Soc., Providence, RI, 1996, pp.~285--293.

\bibitem{DEMMEL}
{\sc J.~W. Demmel}, {\em Applied numerical linear algebra}, Society for
  Industrial and Applied Mathematics (SIAM), Philadelphia, PA, 1997.

\bibitem{GRAU}
{\sc A.~A. Grau}, {\em On the reduction of number range in the use of the
  {G}raeffe process}, J. Assoc. Comput. Mach., 10 (1963), pp.~538--544.

\bibitem{HENRICI}
{\sc P.~Henrici}, {\em Applied and computational complex analysis. {V}ol. 3},
  Pure and Applied Mathematics, John Wiley \& Sons Inc., New York, 1986.
\newblock Discrete Fourier analysis---Cauchy integrals---construction of
  conformal maps---univalent functions, A Wiley-Interscience Publication.

\bibitem{HIGHAM}
{\sc N.~J. Higham}, {\em Accuracy and stability of numerical algorithms},
  Society for Industrial and Applied Mathematics (SIAM), Philadelphia, PA,
  1996.

\bibitem{HOUSEHOLDER}
{\sc A.~S. Householder}, {\em Dandelin, {L}oba\v cevski\u\i, or {G}raeffe?},
  Amer. Math. Monthly, 66 (1959), pp.~464--466.

\bibitem{JS}
{\sc P.~Jana and B.~Sinha}, {\em Fast parallel algorithms for {G}raeffe' s root
  squaring}, Comput. Math. Appl., 35 (1998), pp.~71--80.

\bibitem{TOMS493}
{\sc M.~A. Jenkins}, {\em Algorithm 493 -- {Z}eros of a real polynomial
  {[C2]}}, ACM Transactions on Mathematical Software, 1 (1975), pp.~178--189.

\bibitem{TOMS419}
{\sc M.~A. Jenkins and J.~F. Traub}, {\em Algorithm 419 -- {Z}eros of a complex
  polynomial {[C2]}}, Comm. of the ACM, 15 (1972), pp.~97--99.

\bibitem{KA}
{\sc L.~V. Kantorovich and G.~P. Akilov}, {\em Functional analysis}, Pergamon
  Press, Oxford, second~ed., 1982.
\newblock Translated from the Russian by Howard L. Silcock.

\bibitem{KIRRINIS}
{\sc P.~Kirrinis}, {\em Partial fraction decomposition in {$\mathbb C(z)$} and
  simultaneous newton iteration for factorization in {$\mathbb C[z]$}.}
\newblock Preprint, Bonn, (1995).

\bibitem{KOSTLAN}
{\sc E.~Kostlan}, {\em Random polynomials and the statistical fundamental
  theorem of algebra}, 1987.
\newblock Preprint, MSRI, 1987.

\bibitem{LANG}
{\sc S.~Lang}, {\em Differentiable Manifolds}, Addison-Wesley Pub. Co.,
  Reading, Mass., 1972.

\bibitem{MACKAY}
{\sc R.~S. MacKay}, {\em Renormalisation in area-preserving maps}, vol.~6 of
  Advanced Series in Nonlinear Dynamics, World Scientific Publishing Co. Inc.,
  River Edge, NJ, 1993.

\bibitem{TCS}
{\sc G.~Malajovich}, {\em On generalized {N}ewton algorithms: quadratic
  convergence, path-following and error analysis}, Theoret. Comput. Sci., 133
  (1994), pp.~65--84.
\newblock Selected papers of the Workshop on Continuous Algorithms and
  Complexity (Barcelona, 1993).

\bibitem{MZ96}
{\sc G.~Malajovich and J.~P. Zubelli}, {\em A fast and stable algorithm for
  splitting polynomials}, Comput. Math. Appl., 33 (1997), pp.~1--23.

\bibitem{MZ97}
{\sc G.~Malajovich and J.~P. Zubelli}, {\em On the geometry of {G}raeffe
  iteration}, 1997.
\newblock Informes de Matem\'atica S\'erie B-118, IMPA.

\bibitem{MCMULLEN}
{\sc C.~T. McMullen}, {\em Complex dynamics and renormalization}, vol.~135 of
  Annals of Mathematics Studies, Princeton University Press, Princeton, NJ,
  1994.

\bibitem{NR}
{\sc C.~A. Neff and J.~H. Reif}, {\em An efficient algorithm for the complex
  roots problem}, J. Complexity, 12 (1996), pp.~81--115.

\bibitem{OSTROWSKI}
{\sc A.~Ostrowski}, {\em Recherches sur la m\'ethode de {G}raeffe et les
  z\'eros des polynomes et des s\'eries de {L}aurent}, Acta Math., 72 (1940),
  pp.~99--155.

\bibitem{OSTROWSKI2}
\leavevmode\vrule height 2pt depth -1.6pt width 23pt, {\em Recherches sur la
  m\'ethode de {G}raeffe et les z\'eros des polynomes et des s\'eries de
  {L}aurent. {C}hapitres {I}{I}{I} et {I}{V}}, Acta Math., 72 (1940),
  pp.~157--257.

\bibitem{PAN96}
{\sc V.~Y. Pan}, {\em Optimal and nearly optimal algorithms for approximating
  polynomial zeros}, Comput. Math. Appl., 31 (1996), pp.~97--138.

\bibitem{PAN}
{\sc V.~Y. Pan}, {\em Solving a polynomial equation: some history and recent
  progress}, SIAM Rev., 39 (1997), pp.~187--220.

\bibitem{PKSAH}
{\sc V.~Y. Pan, M.-H. Kim, A.~Sadikou, X.~Huang, and A.~Zheng}, {\em On
  isolation of real and nearly real zeros of a univariate polynomial and its
  splitting into factors}, J. Complexity, 12 (1996), pp.~572--594.
\newblock Special issue for the Foundations of Computational Mathematics
  Conference (Rio de Janeiro, 1997).

\bibitem{PS}
{\sc F.~P. Preparata and M.~I. Shamos}, {\em Computational geometry}, Texts and
  Monographs in Computer Science, Springer-Verlag, New York, 1985.
\newblock An introduction.

\bibitem{RJ}
{\sc T.~A. Rice and L.~H. Jamieson}, {\em A highly parallel algorithm for root
  extraction}, IEEE Trans. Comput., 38 (1989), pp.~443--449.

\bibitem{SCHONHAGE}
{\sc A.~Sch{\"o}nhage}, {\em Equation solving in terms of computational
  complexity}, in Proceedings of the International Congress of Mathematicians,
  Vol. 1, 2 (Berkeley, Calif., 1986), Providence, RI, 1987, Amer. Math. Soc.,
  pp.~131--153.

\bibitem{BEZI}
{\sc M.~Shub and S.~Smale}, {\em Complexity of {B}\'ezout's theorem. {I}.
  {G}eometric aspects}, J. Amer. Math. Soc., 6 (1993), pp.~459--501.

\bibitem{BEZII}
{\sc M.~Shub and S.~Smale}, {\em Complexity of {B}ezout's theorem. {I}{I}.
  {V}olumes and probabilities}, in Computational algebraic geometry (Nice,
  1992), vol.~109 of Progr. Math., Birkh\"auser Boston, Boston, MA, 1993,
  pp.~267--285.

\bibitem{BEZIII}
{\sc M.~Shub and S.~Smale}, {\em Complexity of {B}ezout's theorem. {I}{I}{I}.
  {C}ondition number and packing}, J. Complexity, 9 (1993), pp.~4--14.
\newblock Festschrift for Joseph F. Traub, Part I.

\bibitem{BEZV}
{\sc M.~Shub and S.~Smale}, {\em Complexity of {B}ezout's theorem. {V}.
  {P}olynomial time}, Theoret. Comput. Sci., 133 (1994), pp.~141--164.
\newblock Selected papers of the Workshop on Continuous Algorithms and
  Complexity (Barcelona, 1993).

\bibitem{BEZIV}
{\sc M.~Shub and S.~Smale}, {\em Complexity of {B}ezout's theorem. {I}{V}.
  {P}robability of success; extensions}, SIAM J. Numer. Anal., 33 (1996),
  pp.~128--148.

\bibitem{SMALE}
{\sc S.~Smale}, {\em Newton's method estimates from data at one point}, in The
  merging of disciplines: new directions in pure, applied, and computational
  mathematics (Laramie, Wyo., 1985), Springer, New York, 1986, pp.~185--196.

\bibitem{VAINBERG}
{\sc M.~M. Vainberg}, {\em Variational methods for the study of nonlinear
  operators}, Holden-Day Inc., San Francisco, Calif., 1964.
\newblock With a chapter on Newton's method by L. V. Kantorovich and G. P.
  Akilov. Translated and supplemented by Amiel Feinstein.

\bibitem{WILKINSON}
{\sc J.~H. Wilkinson}, {\em Rounding errors in algebraic processes},
Prentice Hall, Englewood Cliffs NJ (1963).

\bibitem{PERF}
{\sc J.~H. Wilkinson}, {\em The perfidious polynomial}, in Studies in numerical
  analysis, vol.~24 of MAA Stud. Math., Math. Assoc. America, Washington, DC,
  1984, pp.~1--28.

\end{thebibliography}

\end{document}